\documentclass[11pt]{amsart}
\usepackage[english]{babel}
\usepackage{amscd,amssymb,amsmath}
\usepackage{graphicx}
\newtheorem{theorem}[subsection]{Theorem}
\newtheorem{lemma}[subsection]{Lemma}
\newtheorem{problem}{Problem}
\newtheorem{corollary}[subsection]{Corallary}
\def\abd{\mathop {\rm \leftharpoondown \!\! \rightharpoondown} }
\numberwithin{equation}{section}
\begin{document}

\title [Optimal interpolation formulas]
{Optimal interpolation formulas in the periodic function space of
S.L. Sobolev }

\author{Kh.M. Shadimetov, A.R. Hayotov, N.H. Mamatova}

\address{Kh.M. Shadimetov, A.R. Hayotov, N.H. Mamatova\\ Institute of mathematics and information technologies,
Tashkent, Uzbekistan.} \email {hayotov@mail.ru,
abdullo\_hayotov@mail.ru}

\begin{abstract}
In this paper the problem of construction of lattice optimal
interpolation formulas in the space $\widetilde{L_2^{(m)}} (0,1)$
is considered. Using S.L. Sobolev's method explicit formulas for
the coefficients of lattice optimal interpolation formulas are
given and the norm of the error functional of lattice optimal
interpolation formulas is calculated. Moreover, connection between
optimal interpolation formula in the space $\widetilde{L_2^{(m)}}
(0,1)$ and optimal quadrature formula in this space is shown.
Finally, numerical results are given.

\textbf{MSC:} 41A05, 41A15, 65D30, 65D32.

\textbf{Keywords:} lattice optimal interpolation formula; error
functional; optimal coefficients; S.L. Sobolev space; periodic
function.
\end{abstract}

\maketitle

\section{Introduction. Statement of the Problem}

In order to find an approximate representation of a function
$\varphi(x)$ by elements of a certain finite dimensional space, it
is possible to use values of this function at some finite set of
points $x_k$,\ $k=1,2,...,N$. The corresponding problem is called
\emph{the interpolation problem}, and the points $x_k$ \emph{the
interpolation nodes}.

There are polynomial and spline interpolations. Now the theory of
spline interpolation is fast developing. Many books are devoted to
the theory of splines, for example, Ahlberg et al \cite{Ahlb67},
de Boor \cite{Boor78}, Schumaker \cite{Schum81}, Laurent
\cite{Lor75}, Attea \cite{Attea92}, Stechkin and Subbotin
\cite{Stech76}, Vasilenko \cite{Vas83}, Arcangeli et al
\cite{Arc04}, Ignatov and Pevniy \cite{Ign91}, Korneichuk et al
\cite{Korn93}, in the numerical analysis literature or Wahba
\cite{Wahba90}, Eubank \cite{Eub88}, Green and Silverman
\cite{GrSi94}, Berlinet and Thomas-Agnan \cite{BerAgnan04} in the
statistical one.

If the exact values $\varphi(x_k)$ of an unknown smooth function
$\varphi(x)$ at the set of points $\{x_k,\ k=1,2,...,N\}$  in an
interval $[a,b]$ are known, it is usual to approximate $\varphi$
by minimizing
\begin{equation}
\label{eq.1} \int_a^b(g^{(m)}(x))^2dx
\end{equation}
in the set of interpolating functions (i.e. $g(x_k)=\varphi(x_k)$,
$k=1,2,...,N$) of the Sobolev space $L_2^{(m)}(a,b)$. It turns out
that the solution is a natural polynomial spline of order $2m$
with knots $x_1,...,x_N$ called the interpolating $D^m$ spline for
the points $(x_k,\varphi(x_k))$. In non periodic case first this
problem was investigated by Holladay \cite{Hol57} for $m=2$ and
the result of Holladay was generalized by de Boor \cite{deBoor63}
for any $m$. In the Sobolev space $\widetilde{L_2^{(m)} }$ of
periodic functions the minimization problem of integrals of type
(\ref{eq.1}) were investigated by I.J.Schoenberg \cite{Schoen64}
and M.Golomb \cite{Golomb68}.

In the present paper we deal with optimal interpolation formulas.
Now we give the statement of the problem of optimal interpolation
formulas following by S.L.Sobolev \cite{Sob61a}. It should be
noted that connection between the minimization problem of
integrals of type (\ref{eq.1}) and the problem of optimal formulas
was shown in \cite{Sob61b}.

First we recall the definition of Sobolev space
$\widetilde{L_2^{(m)} }(0,1)$ of periodic functions (see
\cite{Sob74,SobVas}).

Suppose in the set $\mathbb{R}$ of real numbers the function
$\varphi (x)$ has local summable derivatives up to order $m$, and
also for the interval $[0,1]$ the integral $\int_0^1 {\left(
{\varphi ^{(m)} (x)} \right)^2 dx}$ is bounded. Assume that the
function $\varphi (x)$ is periodic: $\varphi (x + \gamma ) =
\varphi (x)$,   $x \in \mathbb{R},\,\,\,\gamma \in \mathbb{Z} $,
where $\mathbb{Z}$ is the set of integer numbers.

Every element of the space $\widetilde{L_2^{(m)} }(0,1)$ is a
class of functions which differ from each other by constant term.
The norm of functions in the space $\widetilde{L_2^{(m)} }(0,1)$
is defined by formula
$$\left\| {\left.
{\varphi} \right|\widetilde{L_2^{(m)} }(0,1)} \right\| = \left(
{\int\limits_0^1 {\left( {\varphi ^{(m)} (x)} \right)^2 dx} }
\right)^{1/2}. $$

Now following \cite{Sob61a} we consider interpolation formula of
the form
\begin{equation} \label{eq.(1.1)}
\varphi (x) \cong P_\varphi  (x) = \sum\limits_{k = 1}^N {C_k (x)
\cdot \varphi (x_k )}
\end{equation}
in the space $\widetilde{L_2^{(m)} }(0,1)$. Here points $x_k \in
[0,1]$ and parameters $C_k (x)$ are respectively called \emph{the
nodes} and \emph{the coefficients} of the interpolation formula
(\ref{eq.(1.1)}).

The difference $\varphi-P_{\varphi}$ is called \emph{the error} of
the interpolation formula (\ref{eq.(1.1)}). The value of this
error at some point $z$ is the linear functional on functions
$\varphi$, i.e.
 $$
(\ell,\varphi)\equiv\varphi(z)-P_{\varphi}(z)=\varphi(z)-\sum\limits_{k=1}^N
C_k(z)\varphi(x_{k})=
 $$
\begin{equation}
 \label{eq.(1.2)}
=\int\limits_{0}^1\left[\left(\delta(x-z)-\sum\limits_{k=1}^N
C_k(z)\delta(x-x_k)\right)*\phi_0(x)\right]\varphi(x)dx,
\end{equation}
where $\delta(x)$ is the Dirac delta-function,
$\phi_0(x)=\sum\limits_{\beta}\delta(x-\beta)$, here $\beta$ takes
all integer values and
\begin{equation}\label{eq.(1.3)}
\ell=\left(\delta(x-z)-\sum\limits_{k=1}^NC_k(z)\delta(x-x_k)\right)*\phi_0(x)
\end{equation}
is \emph{the error functional} of interpolation formula
(\ref{eq.(1.1)}) and belongs to the space
$\widetilde{L_2^{(m)*}}(0,1)$. The space
$\widetilde{L_2^{(m)*}}(0,1)$ is the conjugate space to the space
$\widetilde{L_2^{(m)}}(0,1)$ and consists of all periodic
functionals (\ref{eq.(1.3)}) which are orthogonal to unity, i.e.
\begin{equation}\label{eq.(1.4)}
(\ell,1)=0.
\end{equation}

By Cauchy-Schwartz inequality
$$
|(\ell,\varphi)|\leq \|\varphi|\widetilde{L_2^{(m)}}(0,1)\|\cdot
\|\ell|\widetilde{L_2^{(m)*}}(0,1)\|
$$
the error (\ref{eq.(1.2)}) of formula (\ref{eq.(1.1)}) is
estimated with the help of the norm
$$\left\| {\ell|\widetilde{L_2^{(m)*}}(0,1) } \right\| = \mathop {\sup
}\limits_{\left\| {\varphi|\widetilde{L_2^{(m)}}(0,1)} \right\| =
1} \left| {\left( {\ell,\varphi} \right)} \right|
$$
of the error functional (\ref{eq.(1.3)}). Consequently, estimation
of the error of interpolation formula (\ref{eq.(1.1)}) on
functions of the space $\widetilde{L_2^{(m)}}(0,1)$ is reduced to
finding the norm of the error functional $\ell$ in the conjugate
space $\widetilde{L_2^{(m)*}}(0,1)$.

Therefore from here we get the first problem.

\begin{problem} Find the norm of the error functional $\ell$ of
interpolation formula (\ref{eq.(1.1)}) in the space
$\widetilde{L_2^{(m)*}}(0,1)$.
\end{problem}

Obviously the norm of the error functional $\ell (x)$ depends on
the coefficients $C_k(z)$ and the nodes $x_{k}$. \emph{By optimal
interpolation formula} we call such formula which the error
functional in given number $N$ of the nodes has the minimum norm
in the space $\widetilde{L_2^{(m)*}}(0,1)$. If the nodes $x_{k}$
are the points of a lattice, i.e. are located at points of the
form $x_{k}=hk$ then such interpolation formula is called
\emph{the lattice interpolation formula}. Here $h$ is small
parameter and is called \emph{the step of the lattice}.

The main goal of the present paper is to construct the lattice
optimal interpolation formula in the space
$\widetilde{L_2^{(m)}}(0,1)$ for the nodes $x_{k}=hk$, i.e., to
find the coefficients $C_k(z)$ satisfying the following equality
 \begin{equation}
 \label{eq.(1.5)}
\left\| {\mathop \ell \limits^ \circ
(x)|\widetilde{L_2^{(m)*}}(0,1)} \right\| = \mathop {\inf
}\limits_{C_k(z) }  \left\| {\mathop \ell \limits
(x)|\widetilde{L_2^{(m)*}}(0,1)} \right\|.
 \end{equation}

Thus in order to construct the lattice optimal interpolation
formula in the space $\widetilde{L_2^{(m)}}(0,1)$ we need to solve
the next problem.

\begin{problem}
Find the coefficients $C_k(z)$ which satisfy equality
(\ref{eq.(1.5)}) when the nodes are located in the lattice, i.e.,
$x_{k}=hk$.
\end{problem}

In the present paper the lattice optimal interpolation formulas
are constructed in the Sobolev space $\widetilde{L_2^{(m)}}(0,1)$
of periodic functions. First such problem was stated and
investigated by S.L. Sobolev in \cite{Sob61a}, where the extremal
function of the interpolation formula was found in the space
$W_2^{(m)}$.

This paper is organized as follows. In Section 2 we determine the
extremal function which corresponds to the error functional
$\ell(x)$ and give a representation of the norm of the error
functional (\ref{eq.(1.3)}). Section 3 is devoted to minimization
of $\left\| \ell \right\|^2 $ with respect to the coefficients
$C_k(z)$. We obtain a system of linear equations  for the
coefficients of the optimal interpolation formula in the space
$\widetilde{L_2^{(m)}}(0,1)$. Moreover, the existence and
uniqueness of the corresponding solution is proved. In Section 4
we consider the problem of construction of lattice optimal
interpolation formulas. In Section 5 we prove some new properties
of the discrete analogue of the differential operator
$d^{2m}/dx^{2m}$. Using these properties explicit formulas for
coefficients of lattice optimal interpolation formulas are found
in Section 6. In Section 7 the norm of the error functional of
lattice optimal interpolation formulas is calculated. Integrating
obtained lattice optimal interpolation formula the optimal
quadrature formula in the space $\widetilde{L_2^{(m)}}(0,1)$ is
obtained in Section 8. Finally, in Section 9 we give formulas for
coefficients which are very useful in practice and we present some
numerical results.

\section{The extremal function and representation\\ of the norm of the error functional}

In this section we solve Problem 1, i.e. we find explicit form of
the norm of $\ell(x)$.

For finding the explicit form of the norm of the error functional
$\ell(x)$ in the space $\widetilde{L_2^{(m)}}(0,1)$ we use concept
of its extremal function which was introduced by S.L.Sobolev
\cite{Sob61a,Sob74}. The function  $u(x)$ from
$\widetilde{L_2^{(m)}}(0,1)$ is called \emph{the extremal
function} for the error functional $\ell (x)$ if the following
equality is fulfilled
$$
\left( {\ell,u} \right) = \left\| {\ell\left| {\widetilde{L_2^{(m)
* } }(0,1)} \right.} \right\| \cdot \left\| {u\left|
{\widetilde{L_2^{(m)} }(0,1)} \right.} \right\|.
$$
The space $\widetilde{L_2^{(m)}}(0,1)$ is Hilbert space and the
inner product in this space is given by formula
$$
\langle {\varphi,\psi} \rangle_m  = \int\limits_0^1 {\varphi
^{(m)} (x) \cdot \psi ^{(m)} (x)dx}.
$$
According to the Riesz theorem any linear continuous functional
$\ell (x)$ in a Hilbert space is represented in the form of a
inner product
\begin{equation} \label{eq.(2.1)}
\left( {\ell,\varphi} \right) = \langle{\psi _\ell,\varphi}
\rangle_m
\end{equation}
for arbitrary function $\varphi(x)$ from
$\widetilde{L_2^{(m)}}(0,1)$. Here $\psi_\ell  (x)$ is a function
from $\widetilde{L_2^{(m)}}(0,1)$ is defined  uniquely by
functional $\ell (x)$ and is the extremal function. Integrating by
parts the expression in the right hand side of (\ref{eq.(2.1)})
and using periodicity of functions $\varphi (x)$ and $\psi _\ell
(x)$ we get the following equality
$$
\left( {\ell,\varphi} \right) = ( - 1)^m \int\limits_0^1 {{{d^{2m}
} \over {dx^{2m} }}\psi _\ell  (x) \cdot \varphi (x)dx}.
$$
Thus the function $\psi _\ell(x)$ is the generalized solution of
the equation
\begin{equation} \label{eq.(2.2)}
{{d^{2m} }\over {dx^{2m} }}\psi _\ell  (x) = ( - 1)^m \ell (x)
\end{equation}
with the boundary conditions
$$
\psi _\ell ^{(\alpha )} (0) = \psi _\ell ^{(\alpha )}
(1),\,\,\,\alpha  = \overline {0,2m - 1}.
$$

For the extremal function the following holds

\begin{theorem}\label{THM2.1}
Explicit expression for the extremal function $\psi _\ell  (x)$ of
the error functional  (\ref{eq.(1.3)}) is defined  by formula
\begin{equation} \label{eq.(2.3)}
\psi _\ell  (x) = ( - 1)^m \left[ {B_{2m} (x - z) - \sum\limits_{k
= 1}^N {C_k (z) \cdot B_{2m} (x - x_k ) + d_0 } } \right],
\end{equation}
where  $B_{2m}(x) = \sum\limits_{\beta  \ne 0} {{{\exp ( - 2\pi
i\beta x)} \over {(2\pi i\beta )^{2m} }}} $ is the Bernoulli
polynomial, $d_0 $ is a constant.
\end{theorem}

\begin{proof} Here we use the following formulas of Fourier transformations given in
\cite{Sob74}
    $$
F[\varphi (x)] = \int_{ - \infty }^\infty  {\exp (2\pi ipx) \cdot
\varphi (x)dx},\ \ F^{ - 1} [\varphi (p)] = \int_{ - \infty
}^\infty  {\exp ( - 2\pi ixp) \cdot \varphi (p)dp}.
$$
Convolution of two functions is defined by formula
$$
f(x) * g(x) = \int\limits_{ - \infty }^\infty  {f(x - y) \cdot
g(y)dy}.
$$

Applying to both sides of equation (\ref{eq.(2.2)}) the Fourier
transformation and known formulas\linebreak
$F[\delta(x-z)]=e^{2\pi i pz}$, $F[\phi_0(x)]=\phi_0(p)$ (see
\cite{Sob74}) we get
\begin{equation}\label{eq.(2.4)}
(2\pi ip)^{2m} F[\psi _\ell  (x)] = ( - 1)^m \left[ {\left( {\exp
(2\pi ipz) - \sum\limits_{k = 1}^N {C_k (z)\exp (2\pi ipx_k )} }
\right)\phi _0 (p)} \right].
\end{equation}
By virtue of (\ref{eq.(1.4)}) the right hand side of
(\ref{eq.(2.4)}) is zero at the origin. Therefore we can divide
both sides of (\ref{eq.(2.4)}) by $(2\pi ip)^{2m}$. The function
$F[\psi_\ell (x)]$ is defined from equation (\ref{eq.(2.4)}) up to
the following expression
$$
(-1)^md_0 \delta (p) +
\sum\limits_{\alpha = 1}^{2m - 1} {d_\alpha D^\alpha \delta (p)}.
$$
But as known a periodic solution of the homogenous equation
corresponding to equation (\ref{eq.(2.2)}) is a constant term then
all terms except $(-1)^md_0 \delta(p)$ should be omitted. Thus
from (\ref{eq.(2.4)}) we get
\begin{equation} \label{eq.(2.5)}
F[\psi _\ell  (x)] = ( - 1)^m d_0 \delta (p) + {{\exp (2\pi ipz)\
\phi _0 (p)} \over {(2\pi p)^{2m} }} - {{\sum\limits_{k = 1}^N
{C_k (z)\exp (2\pi ipx_k )\ \phi _0 (p)} } \over {(2\pi p)^{2m}
}}.
\end{equation}
Changing the function $\phi _0(p)$ by series of $\delta $
functions and applying the inverse Fourier transformation to both
sides of (\ref{eq.(2.5)}) we get (\ref{eq.(2.3)}). Theorem
\ref{THM2.1} is proved.
\end{proof}

Now we obtain representation for the norm of the error functional
$\ell$. Since the space $\widetilde{L_2^{(m)} }(0,1)$ is the
Hilbert space then by the Riesz theorem we have
\begin{equation}\label{eq.(2.6)}
 \left( {\ell,\psi_\ell} \right) = \left\| {\ell\left|
{\widetilde{L_2^{(m)* } }(0,1)} \right.} \right\| \cdot \left\|
{\psi _\ell \left| {\widetilde{L_2^{(m)} }(0,1)} \right.} \right\|
= \left\| {\ell \left| {\left| {\widetilde{L_2^{(m) * } }(0,1)}
\right.} \right.} \right\|^2 .
\end{equation}
Using formulas (\ref{eq.(1.3)}), (\ref{eq.(2.3)}),
(\ref{eq.(2.6)}) we get
\begin{eqnarray*}
\left\| {\ell(x) \left| {\widetilde{L_2^{(m) * } }(0,1)} \right.}
\right\|^2  &=& \int\limits_0^1 {\ell (x)\psi _\ell (x)dx = } ( -
1)^m \int\limits_0^1 {\left( {\delta (x - z) - \sum\limits_{k =
1}^N {C_k (z) \cdot \delta (x - x_k )} } \right)
* \phi _0 (x)}  \\
&&\times \left( {B_{2m} (x - z) - \sum\limits_{k = 1}^N {C_k (z)
\cdot B_{2m} (x - x_k ) + d_0 } } \right)\\
&=& ( - 1)^m \int\limits_0^1 {\left( {\delta (x - z) * \phi _0 (x)
- \sum\limits_{k = 1}^N {C_k (z)\delta (x - x_k ) * \phi _0 (x)} }
\right) }\\
&&\times \left( {B_{2m} (x - z) - \sum\limits_{k = 1}^N {C_k (z)
\cdot B_{2m} (x - x_k )} } \right)dx + ( - 1)^m d_0 \left( {\ell
(x),1} \right).
\end{eqnarray*}
Hence, using the condition (\ref{eq.(1.4)}) we have
$$
\left\| {\ell(x) \left| {\widetilde{L_2^{(m) * } }(0,1)} \right.}
\right\|^2  = ( - 1)^m \int\limits_0^1 {\left[ {\sum\limits_\beta
{\int\limits_{ - \infty }^\infty  {\delta (x - z - y)\delta (y -
\beta )dy}  - } } \right.}
$$
$$
\left. { - \sum\limits_{k = 1}^N {C_k (z)\sum\limits_\beta
{\int\limits_{ - \infty }^\infty  {\delta (x - x_k  - y)\delta (y
- \beta )dy} } } } \right]\left[ {B_{2m} (x - z) - \sum\limits_{k
= 1}^N {C_k (z) \cdot B_{2m} (x - x_k )} } \right]dx.
$$
Taking into account definition of Dirac's delta-function we obtain
$$
\left\| {\ell(x)\left| {\widetilde{L_2^{(m) * } }(0,1)} \right.}
\right\|^2  = ( - 1)^m \int\limits_0^1 {\left[ {\sum\limits_\beta
{\delta (x - z - \beta )}  - \sum\limits_{k = 1}^N {C_k (z)\delta
(x - x_k  - \beta )} } \right] \times }
$$
$$
\times \left( {B_{2m} (x - z) - \sum\limits_{k = 1}^N {C_k (z)
\cdot B_{2m} (x - x_k )} } \right)dx =
$$
$$
 = ( - 1)^m \left[ {\int\limits_0^1 {\sum\limits_\beta  {\delta
(x - z - \beta ) \cdot B_{2m} (x - z)dx}  - } } \right.
 \sum\limits_{k = 1}^N {C_k (z)\int\limits_0^1
{\sum\limits_\beta  {\delta (x - z - \beta ) \cdot B_{2m} (x - x_k
)dx} }  - }
$$
$$
- \sum\limits_{k = 1}^N {C_k (z)\int\limits_0^1
{\sum\limits_\beta  {\delta (x - x_k  - \beta ) \cdot B_{2m} (x -
z)dx} }  + }
$$
$$
+ \left. {\sum\limits_{k = 1}^N {\sum\limits_{k' = 1}^N {C_k
(z)C_{k'} (z)} \sum\limits_\beta ^{} {\int\limits_0^1 {\delta (x -
x_k  - \beta )B_{2m} (x - x_{k'} )dx} } } } \right].
$$
Hence with the help of the characteristic function $\chi _{\left[
{0,1} \right]} (x)$ of the interval $[0,1]$  square of the norm of
the error functional (\ref{eq.(1.3)}) we reduce to the form
\begin{eqnarray*}
\left\| {\ell(x) \left| {\widetilde{L_2^{(m) * } }(0,1)} \right.}
\right\|^2& =& ( - 1)^m \left[ {\sum\limits_\beta ^{}
{\int\limits_{ - \infty }^\infty  {\chi _{\left[ {0,1} \right]}
(x)\delta (x - z - \beta )\,B_{2m} (x - z)dx  } } } \right.\\
&& - \sum\limits_{k = 1}^N {C_k (z)\sum\limits_\beta ^{}
{\int\limits_{ - \infty }^\infty  {\chi _{[0,1]} (x)\delta (x - z
- \beta )\,B_{2m} (x - x_k )dx  } } }\\
&& - \sum\limits_{k = 1}^N {C_k (z)\sum\limits_\beta ^{}
{\int\limits_{ - \infty }^\infty  {\chi _{[0,1]} (x)\delta (x -
x_k  - \beta )\,B_{2m} (x - z)dx  } } }\\
&& + \left. {\sum\limits_{k = 1}^N {\sum\limits_{k' = 1}^N {C_k
(z)C_{k'} (z)\sum\limits_\beta ^{} {\int\limits_{ - \infty
}^\infty  {\chi _{[0,1]} (x)\delta (x - x_k  - \beta )\,B_{2m} (x
- x_k )dx} } } } } \right].
\end{eqnarray*}
From here applying definition of Dirac's delta-function we have
\begin{eqnarray*}
\left\| {\ell(x) \left| {\widetilde{L_2^{(m) * } }(0,1)} \right.}
\right\|^2 & =& ( - 1)^m \left[ {\sum\limits_\beta  {\chi _{[0,1]}
(z + \beta )B_{2m} (\beta )}  } \right.\\
&& - \sum\limits_{k = 1}^N {C_k (z)} \sum\limits_\beta  {\chi
_{[0,1]} (z + \beta )} B_{2m} (z + \beta  - x_k ) \\
&& - \sum\limits_{k = 1}^N {C_k (z)} \sum\limits_\beta  {\chi
_{[0,1]} (x_k  + \beta )} B_{2m} (x_k  + \beta  - z) \\
&& + \left. {\sum\limits_{k = 1}^N {C_k (z)} \sum\limits_{k' =
1}^N {C_{k'} (z)} \sum\limits_\beta  {\chi _{[0,1]} (x_k  + \beta
)} B_{2m} (x_k  + \beta  - x_{k'} )} \right].
\end{eqnarray*}
Taking into account that $x_k  \in [0,1],\,$ $\,z \in [0,1],\,$
$\,\sum_\beta {\chi _{[0,1]} (y + \beta ) = 1} $ (see.
\cite{Sob74}),  $y \in [0,1]$, and 1-periodicity, symmetry of
$B_{2m} (x)$, that is $B_{2m} (x + \gamma ) = B_{2m}
(x),\,\,\,B_{2m} ( - x) = B_{2m} (x)$ from the last equality we
get the following representation of $\|\ell(x)\|^2$
\begin{eqnarray}
 \left\| {\ell(x)\left| {\widetilde{L_2^{(m) * } }(0,1)} \right.}
\right\|^2 & =& ( - 1)^m \left[ {B_{2m} (0) - 2\sum\limits_{k =
1}^N {C_k (z)B_{2m} (z - x_k ) + } } \right.\nonumber \\
&&  + \left. {\sum\limits_{k = 1}^N {C_k (z)\sum\limits_{k' = 1}^N
{C_{k'} (z)} B_{2m} (x_k  - x_{k'} )} } \right]    .
\label{eq.(2.7)}
\end{eqnarray}

Thus Problem 1 is solved.

Further in next sections we solve Problem 2.

\section{Existence and uniqueness of optimal interpolation formula}

For finding the minimum of the norm (\ref{eq.(2.7)}) under the
condition (\ref{eq.(1.4)}) we use the Lagrange methods. For this
we consider the function
$$
\Psi(\mathbf{C}(z),\lambda)=\|\ell\|^2+2(-1)^m\lambda\
(\ell(x),1),
$$
where $\mathbf{C}(z)=(C_1(z),C_2(z),...,C_N(z))$.

Equating the partial derivatives of $\Psi(\mathbf{C}(z),\lambda)$
by $C_k(z),\ (k=1,2,...,N)$ and $\lambda$ to zero we get the
following system of linear equations
\begin{eqnarray}
&&
\sum\limits_{k=1}^NC_k(z)B_{2m}(x_{k'}-x_{k})+\lambda=B_{2m}(z-x_{k'}),
\  \ k'=1,2,...,N,
 \label{eq.(3.1)} \\
&& \sum\limits_{k=1}^NC_k(z)=1.
 \label{eq.(3.2)}
\end{eqnarray}

System (3.1)-(3.2) has a unique solution and this solution gives
the minimum to $\|\ell(x)\|^2$ under the condition
(\ref{eq.(3.2)}).

It should be noted that uniqueness of the solution of such type
systems were also investigated in
\cite{Vas83,Arc04,Ign91,Sob74,SobVas}.

The uniqueness of the solution of system
(\ref{eq.(3.1)})--(\ref{eq.(3.2)})  is proved following
\cite[Chapter I]{SobVas}. For completeness we give it here.

First in (\ref{eq.(2.7)}) we change of variables
$C_k(z)=\bar{C}_k(z)+d_k(z)$ then (\ref{eq.(2.7)}) and system
(\ref{eq.(3.1)})--(\ref{eq.(3.2)}) have the following form
\begin{eqnarray}
\|\ell\|^2&=& (-1)^m\left[\sum\limits_{k=1}^N\bar{C}_k(z) \right.
  \sum\limits_{k'=1}^N\bar{C}_{k'}(z)B_{2m}(x_{k}-x_{k'})-2\sum\limits_{k=1}^N
(\bar{C}_k(z)+d_k(z))B_{2m}(z-x_{k})\nonumber \\
&&+ \sum\limits_{k=1}^N\sum\limits_{k'=1}^N
 \left.(d_k(z)d_{k'}(z)+\bar{C}_k(z)d_{k'}(z)+\bar{C}_{k'}(z)d_k(z))
B_{2m}(x_{k}-x_{k'})+B_{2m}(0)\right], \label{eq.(3.3)}
\end{eqnarray}
\begin{eqnarray}
&&
\sum\limits_{k=1}^N\bar{C}_k(z)B_{2m}(x_{k'}-x_{k})+\lambda=B_{2m}(z-x_{k'})-\sum\limits_{k=1}^Nd_k(z)B_{2m}(x_{k'}-x_{k}),
 \label{eq.(3.4)}  \\
&&  k'=1,2,...,N,  \nonumber \\
&&\sum\limits_{k=1}^N \bar{C}_k(z)=0,\label{eq.(3.5)}
\end{eqnarray}
where $d_k(z)$ is a partial solution of the equation
(\ref{eq.(3.2)}).

Hence we directly get that the minimization of (2.7) under the
condition (\ref{eq.(3.2)}) with respect to $C_k(z)$ is equivalent
to the minimization of expression (\ref{eq.(3.3)}) with respect to
$\overline C _k(z)$ under the condition (\ref{eq.(3.5)}).
Therefore it is sufficient to prove that system
(\ref{eq.(3.4)})--(\ref{eq.(3.5)}) has a unique solution with
respect to unknowns $\overline{\mathbf{C}}(z) = (\overline C _1(z)
,\overline C _2(z),...,\overline C _N(z)),$ $\lambda$ and this
solution gives the conditional minimum for
$\left\|\ell\right\|^2$.

From the theory of conditional extremum it is known the sufficient
condition in which the solution of system
(\ref{eq.(3.4)})--(\ref{eq.(3.5)}) gives the conditional minimum
for $\left\| \ell \right\|^2$ on manifold (\ref{eq.(3.5)}). It
consists of positiveness of the following quadratic form
\begin{equation}
\Phi(\overline{\mathbf{C}}(z))=\sum\limits_{k=1}^N \sum\limits_{k'
=1}^N {{{\partial ^2 \Psi(\overline{\mathbf{C}}(z),\lambda) }
\over {\partial \overline C _k(z)
\partial \overline C _{k'}(z)}}}\overline C_k(z) \overline C
_{k'}(z)
\label{eq.(3.6)}
\end{equation}
on the set of vectors  $\overline{\mathbf{C}}(z)  = (\overline {C}
_1(z),\overline C_2(z),...,\overline C_N(z) )$ under the condition
\begin{equation}
S\overline{\mathbf{C}}  = 0,\label{eq.(3.7)}
\end{equation}
where $S=(1,1,...,1)$ is the $N$ component vector.

In our case this condition is fulfilled, i.e. the following holds

\begin{lemma}\label{L3.1} For any non zero vector
$\overline{\mathbf{C}}(z)=(\overline C_1(z),\overline
C_2(z),...,\overline C_N(z))\in \mathbb{R}^N$ lying in the
subspace (\ref{eq.(3.7)}) the function
$\Phi(\overline{\mathbf{C}}(z))$ is strictly positive.
\end{lemma}

\begin{proof}
From definition of the function
$\Psi(\overline{\mathbf{C}}(z),\lambda)$ and (\ref{eq.(3.6)}) it
follows that
\begin{equation}
\Phi(\overline{\mathbf{C}}(z))=2\sum_{k=1}^N\sum_{k'=1}^N(-1)^mB_{2m}(x_{k}-x_{k'})\overline
C_k(z)\overline C_{k'}(z).\label{eq.(3.8)}
\end{equation}
We consider the functional
$$
\ell_{\overline{\mathbf{C}}(z)}(x)=\sqrt{2}\sum_{k=1}^N\overline
C_k(z)\delta(x-x_{k})*\phi_0(x).
$$
By virtue of condition (\ref{eq.(3.7)}) the functional
$\ell_{\overline{\mathbf{C}}(z)}(x)$ belongs to the space
$\widetilde{L_2^{(m)*}}(0,1)$. Thus this functional has the
extremal function $u_{\overline{\mathbf{C}}(z)}(x)\in
\widetilde{L_2^{(m)}}(0,1)$ which is the solution of the equation
\begin{equation}
\frac{d^{2m}}{dx^{2m}}u_{\overline{\mathbf{C}}(z)}(x)=(-1)^m\ell_{\overline{\mathbf{C}}(z)}(x).
\label{eq.(3.9)}
\end{equation}
As $u_{\overline{\mathbf{C}}(z)}(x)$ we take the following
function
$$
u_{\overline{\mathbf{C}}(z)}(x)=\sqrt{2}\sum_{k=1}^N(-1)^mC_k(z)B_{2m}(x-x_{k}).
$$
Square of the norm of $u_{\overline{\mathbf{C}}(z)}(x)$ coincide
with the function $\Phi(\overline{\mathbf{C}}(z))$ in the space
$\widetilde{L_2^{(m)}}(0,1)$, i.e.
\begin{equation}
\|u_{\overline{\mathbf{C}}(z)}(x)|\widetilde{L_2^{(m)}}(0,1)\|^2=\left(\ell_{\overline{\mathbf{C}}(z)}(x),u_{\overline{\mathbf{C}}(z)}(x)\right)=
2\sum_{k=1}^N\sum_{k'=1}^N(-1)^mB_{2m}(x_{k}-x_{k'})\overline{C}_k(z)\overline{C}_{k'}(z).
\label{eq.(3.10)}
\end{equation}
From here taking into account (\ref{eq.(3.8)}) we conclude that
for non zero $\overline{\mathbf{C}}(z)$ the function
$\Phi(\overline{\mathbf{C}}(z))$ is strictly positive. Lemma
\ref{L3.1} is proved.\end{proof}

If the system (\ref{eq.(3.4)})--(\ref{eq.(3.5)}) has a unique
solution then the system (\ref{eq.(3.1)})--(\ref{eq.(3.2)}) has
also a unique solution.

The following holds

\begin{lemma}\label{L3.2}
The main matrix $Q$ of the system
(\ref{eq.(3.4)})--(\ref{eq.(3.5)}) is non singular.
\end{lemma}

\begin{proof} The homogenous system corresponding to system (\ref{eq.(3.4)})--(\ref{eq.(3.5)}) have the
following matrix form
\begin{equation}
Q\left(%
\begin{array}{c}
  \overline{\mathbf{C}}(z) \\
  \lambda \\
\end{array}%
\right)
=\left(%
\begin{array}{cc}
  B & 1 \\
  S & 0 \\
\end{array}%
\right)
\left(%
\begin{array}{c}
  \overline{\mathbf{C}}(z) \\
  \lambda \\
\end{array}%
\right) =0, \label{eq.(3.11)}
\end{equation}
where $B$ is the matrix with elements
$b_{ij}=B_{2m}(x_{i}-x_{j}),\ \ i=1,2,...,N ,$  $j=1,2,...,N,$ $S$
is the vector $1\times N$ and $S=(1,1,...,1)$.

We show that system (\ref{eq.(3.11)}) has the unique solution
$\overline{\mathbf{C}}(z)=(0,0,...,0),$ $\lambda=0. $

Assume $\overline{\mathbf{C}}(z),$ $\lambda$ is a solution of
(\ref{eq.(3.11)}). Since a solution of equation (\ref{eq.(3.9)})
is determined up to constant term then as
$u_{\overline{\mathbf{C}}(z)}(x)$ we can take the following
function
 $$
u_{\overline{\mathbf{C}}(z)}(x)=(-1)^m \sqrt{2}\left(
\sum\limits_{k'=1}^N
\overline{C}_{k'}(z)B_{2m}(x-x_{k'})+\lambda\right).
 $$
But from (\ref{eq.(3.11)}) it is clear that
$u_{\overline{\mathbf{C}}(z)}(x_{k})=0$. Then for the norm of the
functional $\ell_{\overline{\mathbf{C}}(z)}(x)$ we have
$$
\left\| {\ell_{\overline{\mathbf{C}}(z)}(x)\left|
{\widetilde{L_2^{(m) *}}(0,1)} \right.} \right\|^2 = \left(
{\ell_{\overline{\mathbf{C}}(z)}(x),u_{\overline{\mathbf{C}}(z)}
(x)}\right) = \sum\limits_{k = 1}^N {\overline{\mathbf{C}}_k (z)
\cdot u_{\overline{\mathbf{C}}(z)}(x_{k}) = 0}
$$
on the other hand from (\ref{eq.(3.10)}) we get
$$
\left\| {\ell_{\overline{\mathbf{C}}(z)}(x)\left|
{\widetilde{L_2^{(m) * } }(0,1)} \right.} \right\|^2 = \left(
{\ell_{\overline{\mathbf{C}}(z)}(x),u_{\overline{\mathbf{C}}(z)}
(x)}\right) =
2\sum_{k=1}^N\sum_{k'=1}^N(-1)^mB_{2m}(x_{k}-x_{k'})\overline{C}_k(z)\overline{C}_{k'}(z),
$$
which is possible only when $\overline C_k(z) = 0,$ $k = \overline {1,N} $.\\
Hence by virtue of (\ref{eq.(3.11)}) we obtain $\lambda  = 0.$
Lemma \ref{L3.2} is proved. \end{proof}

From (\ref{eq.(2.7)}) and Lemmas 3.1, 3.2 it follows that in fixed
values of the nodes $x_{k}$ square of the norm of the error
functional $\ell(x)$ being quadratic function of the coefficients
$C_k(z)$ has a unique minimum in some concrete value $C_k(z)=
\mathop {C_k }\limits^ \circ(z)$. Interpolation formulas with the
coefficients $\mathop {C_k }\limits^ \circ(z) \,\,(k  =
\overline{1,N})$, corresponding to this minimum in fixed values of
the nodes $x_{k}$ is called \emph{the optimal interpolation
formulas} and $\mathop {C_k }\limits^ \circ(z) \,\,(k= \overline
{1,N} )$ are called \emph{the optimal coefficients}.

\section{The system for the coefficients of lattice optimal interpolation
formula}

We consider system (\ref{eq.(3.1)})--(\ref{eq.(3.2)}) from Section
3 on one dimensional lattice, i.e., suppose that the nodes
$x_\gamma =h\gamma $, where $h = {1 \over N},\,$ $\,\gamma  = 1,2,
\ldots,N$. Below for convenience we denote $[\gamma]=h\gamma$.
Then such an interpolation formula we call \emph{the lattice
interpolation formula}. Moreover in this case system
(\ref{eq.(3.1)})--(\ref{eq.(3.2)}) takes the following form
\begin{eqnarray}
& & \sum\limits_{\gamma  = 1}^N {\mathop C\limits^ \circ  \left(
{[\gamma ],z} \right)B_{2m} {[\beta  - \gamma ]} + \lambda =
B_{2m}\left( {z - [\beta]} \right),\,\,\,\,\,\,\,}
\beta  = 1,2,...,N, \label{eq.(4.1)}  \\
&&\sum\limits_{\gamma  = 1}^N {\mathop C\limits^ \circ  \left(
{[\gamma ],z} \right)}  = 1. \label{eq.(4.2)}
\end{eqnarray}

Later we use the theory of discrete argument functions.  The
theory of discrete argument functions was investigated in
\cite{Sob74}.

The convolution of two discrete argument functions is given by
formula (cf. \cite{Sob74})
$$
f[\beta] * g[\beta ] = \sum\limits_{\gamma  =  - \infty }^\infty
{f[\gamma ]}  \cdot g[\beta  - \gamma ].
$$

Below for convenience instead of the sum $\sum\nolimits_{\gamma  =
- \infty }^\infty$ we write $\sum\nolimits_{\gamma}$.

Using the discrete characteristic function $\chi _{[0,1]}[\beta ]$
of the interval $[0,1]$ and taking into account the definition of
convolution of two discrete argument functions system
(\ref{eq.(4.1)})--(\ref{eq.(4.2)}) we rewrite in the following
convolution form
\begin{eqnarray}
&& B_{2m} [\beta ] * \left( {\mathop C\limits^ \circ  ([\beta ],z)
\cdot \chi _{[0,1]} [\beta ]} \right) + \lambda  = B_{2m} \left(
{z - [\beta ]} \right),\,\,\,\,\,\,\,\,\,\,[\beta ] \in (0,1],
\label{eq.(4.3)}\\
&& \sum\limits_{\beta  = 1}^N {\mathop C\limits^ \circ  \left(
{[\beta ],z} \right) = 1,}\label{eq.(4.4)}
\end{eqnarray}
where $[\beta ] = h\beta  = {\beta  \over N}$.

Now we have the following problem.

\begin{problem}
Find a discrete function $\mathop C\limits^\circ \left({[\beta
],z} \right)$ and unknown constant $\lambda$ which satisfy the
system (\ref{eq.(4.3)})--(\ref{eq.(4.4)}).
\end{problem}

In the solution of Problem 3 the main role plays some new property
of the discrete analogue $D_h^{(m)}[\beta]$ of the differential
operator $d^{2m}/dx^{2m}$. It should be noted that the properties
of the discrete analogue of the polyharmonic operator $\Delta^m$
were investigated by S.L. Sobolev \cite{Sob74,Sob65}. But here we
need some new properties of the discrete operator
$D_h^{(m)}[\beta]$. The next section is devoted to investigation
of these properties.

\section{Some new properties of the discrete operator  $D_h^{(m)}
[\beta ]$}

In the work \cite{Shad85} was constructed the discrete analogue
$D_h^{(m)}[\beta ]$ of the differential operator ${{d^{2m} } \over
{dx^{2m} }}$. The discrete operator $D_h^{(m)}[\beta ]$ is the
solution of the following difference equation
\begin{equation}
D_h^{(m)}[\beta ]* G_m [\beta ] = \delta [\beta ].
\label{eq.(5.1)}
\end{equation}
Here $G_m [\beta ] = {{\left| {h\beta } \right|^{2m - 1} } \over
{2 \cdot (2m - 1)!}}$,   $\delta [\beta ]$  is equal to 1 when
$[\beta ]$ =0 and is equal to zero when $[\beta ] \ne 0$. The
discrete operator $D_h^{(m)}[\beta ]$ which satisfies equality
(\ref{eq.(5.1)}) has the following form
\begin{equation}
D_h^{(m)} [\beta ] = {{(2m - 1)!} \over {h^{2m} }}\left\{
\begin{array}{ll}
\sum\limits_{k = 1}^{m - 1} {{(1 - q_k )^{2m + 1} \cdot
q_k^{\left| \beta  \right|} } \over {q_k  \cdot E_{2m - 1} (q_k
)}},&\left| \beta  \right| \ge 2,  \cr 1 + \sum\limits_{k = 1}^{m
- 1} {{(1 - q_k )^{2m - 1} } \over {E_{2m - 1} (q_k )}}, & \left|
\beta  \right| = 1,  \cr - 2^{2m - 1}  + \sum\limits_{k = 1}^{m -
1} {{(1 - q_k )^{2m + 1} } \over {q_k \cdot E_{2m - 1} (q_k
)}},&\beta  = 0, \cr
\end{array}
 \right.\label{eq.(5.2)}
\end{equation}
where $\left|{q_k } \right| < 1$ are the roots of the
Euler-Frobenius polynomial $E_{2m - 2} (\lambda )$ of degree
$2m-2$.

The definition of the Euler-Frobenius polynomial is given, for
example, in \cite{SobVas}.

The following holds

\begin{lemma}\label{L5.1}
 The function $D_h^{(m)}[\beta ]$ is the solution
of the equation
\begin{equation}
 hD_h^{(m)} [\beta ] * B_{2m} [\beta ] = \Phi [\beta ] -
 h,\label{eq.(5.3)}
\end{equation}
where $B_{2m} [\beta ] = \sum\limits_{\gamma  \ne 0} {{{\exp
(-2\pi i\gamma h\beta )} \over {(2\pi i\gamma )^{2m} }}} $,
\begin{equation}
 \Phi [\beta ] = \sum\limits_\gamma  {\delta [\beta  - \gamma h^{
- 1} ]} ,\label{eq.(5.4)}
\end{equation}
here $\delta [\beta  - h^{ - 1} \gamma ]$  is equal to 1 when
$\beta = \gamma h^{ - 1} $ and is equal to zero when $\beta \ne
\gamma h^{-1} $.
\end{lemma}

\begin{proof}
We investigate the solution of equation  (\ref{eq.(5.3)}). For
this we use the Fourier  transformation. For convenience first we
go from functions of discrete argument to harrow-shaped functions.
According to the definition of harrow-shaped functions (see
\cite{Sob74,SobVas}) for arbitrary function $\psi [\beta ]$ of
discrete argument its harrow-shaped function has the form
$$
\stackrel{\abd}{\psi}\!\!(x) = \sum\limits_\beta  {h\psi [\beta
]\delta (x - h\beta )}.
$$
For $\Phi [\beta ]$ which defined by formula (\ref{eq.(5.4)}) we
have
\begin{eqnarray*}
\stackrel{\abd}{\Phi}\!\!(x) &=& \sum\limits_\beta  {h\Phi [\beta
]\delta (x - h\beta )}  = \sum\limits_\beta {h\sum\limits_\gamma
{\delta [\beta - h^{ - 1} \gamma ]} \delta (x - h\beta )}\\
& =& \sum\limits_\gamma  {h\delta (x - \gamma )}  = h\phi _0 (x),
\end{eqnarray*}
where $\phi _0 (x) = \sum\limits_\gamma  {\delta (x - \gamma
)}$.\\
Taking into account the last equality and the identity $\delta
(hx) = h^{ - 1} \delta (x)$ for equation (\ref{eq.(5.3)}) in the
class of harrow-shaped functions we have
\begin{equation}
h\stackrel{\abd}{D}_h^{(m)}\!\!(x)
*\stackrel{\abd}{B}_{2m}\!\!(x)= h\phi _0 (x) - h\phi _0 (xh^{ -
1} ).\label{eq.(5.6)}
\end{equation}
The Fourier transformation of the functions $\phi _0 (x)$ and
$\phi _0 (xh^{ - 1} )$ are respectively given  by formulas
\begin{equation}
F[\phi _0 (x)] = \int\limits_{ - \infty }^\infty  {e^{2\pi ipx} }
\sum\limits_\beta  {\delta (x - \beta )dx = } \sum\limits_\beta
{e^{2\pi ip\beta } },\label{eq.(5.7)}
\end{equation}
\begin{equation}
F[\phi _0 (xh^{ - 1} )] = \int\limits_{-\infty}^{\infty} {e^{2\pi
ipx} } \sum\limits_\beta {\delta (xh^{ - 1}  - \beta )dx = }
h\sum\limits_\beta  {e^{2\pi iph\beta } } .\label{eq.(5.8)}
\end{equation}
Hence using the following known formula  (see \cite{Sob74})
\begin{equation}
h\sum\limits_\beta  {e^{2\pi iph\beta } }  = \sum\limits_\beta
{\delta (p - h^{ - 1} \beta )}       \label{eq.(5.9)}
\end{equation}
equalities (\ref{eq.(5.7)}), (\ref{eq.(5.8)}) we reduce to the
form
\begin{eqnarray}
F[\phi _0 (x)] &=& \phi _0 (p),\label{eq.(5.10)}\\
F[\phi _0 (xh^{ - 1} )] &=& \sum\limits_\beta  {\delta (p - h^{ -
1} \beta )dx = } h\phi _0 (hp).\label{eq.(5.11)}
\end{eqnarray}
Now we calculate the Fourier transformation of the function
$$\stackrel{\abd}{B}_{2m} (x) = ( - 1)^m \sum\limits_\beta  h
\sum\limits_{\gamma  \ne 0} {{{\exp ( - 2\pi i\gamma h\beta
)\delta (x - h\beta )} \over {(2\pi \gamma )^{2m} }}}.
$$
By definition of the Fourier transformation we have
\begin{eqnarray*}
F[\stackrel{\abd}{B}_{2m}\!\!(x)] &= &\int\limits_{- \infty
}^\infty {e^{2\pi ipx} }\stackrel{\abd}{B}_{2m}\!\! (x)dx = {{( -
1)^m } \over {(2\pi )^{2m} }}\sum\limits_\beta  h
\sum\limits_{\gamma \ne 0} {{{\exp ( - 2\pi i\gamma h\beta )\exp
(2\pi iph\beta )} \over {\gamma ^{2m} }}}\\
& =& {{( - 1)^m h} \over {(2\pi )^{2m} }}\sum\limits_\beta
{\sum\limits_{\gamma  \ne 0} {{{\exp (2\pi ih\beta (p - \gamma ))}
\over {\gamma ^{2m} }}} }.
\end{eqnarray*}
By virtue of equality (\ref{eq.(5.9)}) we obtain
$$
h\sum\limits_\beta  {\exp (2\pi ih\beta (p - \gamma ))}  =
\sum\limits_\beta  {\delta (p - \gamma  - h^{ - 1} \beta )}.
$$
Therefore
$$
F[\stackrel{\abd}{B}_{2m}\!\!(x)] = {{( - 1)^m } \over {(2\pi
)^{2m} }}\sum\limits_\beta  {\sum\limits_{\gamma  \ne 0} {{{\delta
(p - \gamma  - h^{ - 1} \beta )} \over {\gamma ^{2m} }}} }.
$$
Hence, setting  $\gamma + h^{ - 1} \beta = k,\,\,\,\gamma  = k -
h^{ - 1} \beta $, we get
\begin{equation}
F[\stackrel{\abd}{B}_{2m} (x)] = {{( - 1)^m } \over {(2\pi )^{2m}
}}\sum\limits_\beta {\sum\limits_{\scriptstyle k \atop
 \scriptstyle kh \notin \mathbb{Z}}
{{{\delta (p - k)} \over {(k - h^{ - 1} \beta )^{2m} }}} .}
\label{eq.(5.12)}
\end{equation}
Now applying the Fourier transformation to both sides of equation
(\ref{eq.(5.6)}),  using (\ref{eq.(5.10)}), (\ref{eq.(5.11)}),
(\ref{eq.(5.12)}) and obtained equality dividing by $h$ we get the
following equation
\begin{equation}
 F[\stackrel{\abd}{D}_h^{(m)}\!\!(x)] \cdot {{( - 1)^m } \over
{(2\pi )^{2m} }}\sum\limits_\beta {\sum\limits_{\scriptstyle
\,\,\,\,\,\gamma  \hfill \atop \scriptstyle h\gamma  \notin
\mathbb{Z}} {{{\delta (p - \gamma )} \over {(\gamma  - h^{ - 1}
\beta )^{2m} }}} }  = \sum\limits_{\scriptstyle \,\,\,\,\,\gamma
\hfill \atop \scriptstyle h\gamma  \notin \mathbb{Z}}  {\delta (p
- \gamma )}.\label{eq.(5.13)}
\end{equation}
It is known that
\begin{equation}
\sum\limits_\beta  {\sum\limits_{\scriptstyle \,\,\,\,\gamma
\hfill \atop \scriptstyle \gamma h \notin \mathbb{Z}} {{{\delta (p
- \gamma )} \over {(\gamma  - h^{ - 1} \beta )^{2m} }}} }  =
\sum\limits_\beta  {\sum\limits_{\scriptstyle \,\,\,\,\gamma
\hfill \atop
  \scriptstyle \gamma h \notin \mathbb{Z}}
{{{\delta (p - \gamma )} \over {(p - h^{ - 1} \beta )^{2m} }}} }
.\label{eq.(5.14)}
\end{equation}
By virtue of equality (\ref{eq.(5.14)}) equation (\ref{eq.(5.13)})
is equivalent to the following equation
\begin{equation}
F[\stackrel{\abd}{D}_h^{(m)}\!\!(x)]\left({\Gamma_h ^{(m)}(p)}
\right)^{ - 1}  = 1,\,\,\,\,p \ne h^{ - 1}
\gamma,\label{eq.(5.15)}
\end{equation}
where
\begin{equation}
\Gamma_h^{(m)} (p) = \left[ {{{( - 1)^m } \over {(2\pi )^{2m}
}}\sum\limits_\gamma  {{1 \over {(p - h^{ - 1} \gamma )^{2m} }}} }
\right]^{ - 1}. \label{eq.(5.16)}
\end{equation}
The function  $\Gamma_h^{(m)}(p)$ is periodic with respect to $p$
with period $h^{ - 1} $, real and analytic in all $p \ne h^{ - 1}
\gamma ,\,\,\gamma  \in \mathbb{Z}$.

From equality (\ref{eq.(5.15)}) we have
\begin{equation}
F[\stackrel{\abd}{D}_h^{(m)}\!\!(x)]=\Gamma_h^{(m)}(p).
\label{eq.(5.17)}
\end{equation}

Applying to both sides of equality (\ref{eq.(5.17)}) the inverse
Fourier transformation and going from harrow-shaped functions to
functions of a discrete argument we find (\ref{eq.(5.2)}) (see
\cite{Shad85}). This completes the proof of Lemma \ref{L5.1}.
\end{proof}

\begin{lemma}\label{L5.2}
For the convolution of discrete functions $D_h^{(m)}[\beta ]$ and
$\exp(2\pi i\sigma h\beta)$  the following holds
\begin{eqnarray}
&& D_h^{(m)}[\beta ]*\exp (2\pi i\sigma h\beta ) = ( - 1)^m (2\pi
)^{2m} \exp (2\pi i\sigma h\beta ) \cdot \left[
{\sum\limits_\gamma  {{{h^{2m} } \over {(\gamma  - \sigma h)^{2m}
}}} } \right]^{ - 1} ,\,\,\,\,\sigma h\not\in
\mathbb{Z},\label{eq.(5.18)}\\
&& D_h^{(m)} [\beta ]* \exp (2\pi i\sigma h\beta ) =
0,\,\,\,\,\sigma h \in \mathbb{Z}.\label{eq.(5.19)}
\end{eqnarray}
where $D_h^{(m)}[\beta ]$ is defined by equality (\ref{eq.(5.2)}).
\end{lemma}

\begin{proof}
The convolution of two harrow-shaped functions
$\stackrel{\abd}{D}_h^{(m)}\!\!(x)$ and
$\stackrel{\abd}{\exp}\!\!(2\pi i\sigma x)$ we denote by
$\stackrel{\abd}{T}\!\!(x)$, i.e
$$
\stackrel{\abd}{T}\!\!(x)=\stackrel{\abd}{D}_h^{(m)}\!\!(x)*\stackrel{\abd}{\exp}\!\!(2\pi
i\sigma x).
$$
Using formula (\ref{eq.(5.17)}) we calculate the Fourier
transformation of the function $\stackrel{\abd}{T}\!\!(x)$
\begin{eqnarray}
F[\stackrel{\abd}{T}(x)] &=& F[\stackrel{\abd}{D}_h^{(m)}\!\!(x) *
\stackrel{\abd}{\exp}\!\!(2\pi i\sigma x)] \nonumber \\
&=& F[\stackrel{\abd}{D}_h^{(m)}\!\!(x)] \cdot
F[\stackrel{\abd}{\exp}\!\!(2\pi
i\sigma x)] \nonumber \\
&=&\Gamma_h^{(m)}(p) \cdot F[\stackrel{\abd}{\exp}\!\!(2\pi
i\sigma x)]. \label{eq.(5.20)}
\end{eqnarray}
By definition of harrow-shaped functions and the Fourier
transformation of $\delta (x - h\beta )$ we have
\begin{eqnarray*}
F[\stackrel{\abd}{\exp}\!\!(2\pi i\sigma x)] &=&
F[\sum\limits_\beta
{h\exp (2\pi i\sigma h\beta )\delta (x - h\beta )} ] \\
&=& \sum\limits_\beta  {h\exp (2\pi i\sigma h\beta )} \exp (2\pi
iph\beta ) \\
&=& \sum\limits_\beta{h\exp (2\pi ih\beta (\sigma+ p))}.
\end{eqnarray*}
Hence taking into account (\ref{eq.(5.9)}) we find
\begin{equation}
 F[\stackrel{\abd}{\exp}(2\pi i\sigma x)]=
\sum\limits_\beta  {\delta (\sigma  + p - h^{ - 1} \beta
)}.\label{eq.(5.21)}
\end{equation}
Taking into account formulas (\ref{eq.(5.16)}), (\ref{eq.(5.21)})
from (\ref{eq.(5.20)}) we get
\begin{eqnarray*}
F[\stackrel{\abd}{T}(x)]&=&\left[ {{{( - 1)^m } \over {(2\pi
)^{2m} }}\sum\limits_\gamma  {{1 \over {(p - h^{ - 1} \gamma
)^{2m} }}} } \right]^{ - 1}  \cdot \sum\limits_\beta {\delta
(\sigma  + p - h^{ - 1} \beta )}  \\
&=& \sum\limits_\beta  {\delta (\sigma  + p - h^{ - 1} \beta )}
\left[ {{{( - 1)^m } \over {(2\pi )^{2m} }}\sum\limits_\gamma  {{1
\over {(h^{ - 1} \beta  - \sigma  - h^{ - 1} \gamma )^{2m} }}} }
\right]^{ - 1}  \\
& =& \sum\limits_\beta  {\delta (\sigma  + p - h^{ - 1} \beta )}
\cdot \left[ {{{( - 1)^m } \over {(2\pi )^{2m}
}}\sum\limits_\gamma  {{1 \over {\left( {h^{ - 1} (\beta  - \gamma
) - \sigma } \right)^{2m} }}} } \right]^{ - 1} ,\,\,\,\beta  -
\gamma  \ne \sigma h.\\
\end{eqnarray*}
From here, since $\beta $ and $\gamma $ take all integer values,
we obtain
\begin{equation}
F[\stackrel{\abd}{T}\!\!(x)]=\sum\limits_\beta {\delta (\sigma + p
- h^{ - 1} \beta )}  \cdot \left[ {{{( - 1)^m } \over {(2\pi
)^{2m} }}\sum\limits_\gamma  {{{h^{2m} } \over {\left( {\gamma  -
\sigma h} \right)^{2m} }}} } \right]^{ - 1} ,\,\,\gamma  \ne
\sigma h, \label{eq.(5.22)}
\end{equation}
i.e. $\sigma h \notin \mathbb{Z}.$

Using the identity (\ref{eq.(5.9)}), we calculate the inverse
Fourier transformation of the function\linebreak
$\sum\limits_\beta {\delta (\sigma  + p - h^{ - 1} \beta )}$, i.e.
\begin{eqnarray}
F^{-1}[\sum\limits_\beta  {\delta (\sigma  + p - h^{ - 1} \beta )}
]& = &\sum\limits_\beta  {\int\limits_{ - \infty }^\infty  {e^{ -
2\pi ixp} } \delta (\sigma  + p - h^{ - 1} \beta )} dp =
\sum\limits_\beta  {\exp \left( { - 2\pi ix(h^{ - 1} \beta  -
\sigma )} \right)} \nonumber \\
 &=& e^{2\pi ix\sigma }
\sum\limits_\beta {e^{ - 2\pi ixh^{ - 1} \beta } }  = e^{2\pi
ix\sigma } \sum\limits_\beta {h\delta (x - h\beta )}\nonumber \\
& =& \sum\limits_\beta  {he^{2\pi i\sigma h\beta } } \delta (x -
h\beta ) =\stackrel{\abd}{\exp}\!\!(2\pi i\sigma
x).\label{eq.(5.23)}
\end{eqnarray}
Now applying the inverse Fourier transformation to both sides of
equality (\ref{eq.(5.22)}) and taking into account
(\ref{eq.(5.23)}) we have
$$
\stackrel{\abd}{T}(x)=\stackrel{\abd}{\exp}\!\!(2\pi i\sigma
x)\cdot\left[ {{{( - 1)^m } \over {(2\pi )^{2m}
}}\sum\limits_\gamma  {{{h^{2m} } \over {\left( {\gamma  - \sigma
h} \right)^{2m} }}} } \right]^{ - 1} ,\,\sigma h \notin
\mathbb{Z}.
$$
Hence, going from harrow-shaped functions to functions of discrete
argument, we obtain
$$
T[\beta ] = \exp (2\pi i\sigma h\beta ) \cdot \left[ {{{( - 1)^m }
\over {(2\pi )^{2m} }}\sum\limits_\gamma {{{h^{2m} } \over {\left(
{\gamma  - \sigma h} \right)^{2m} }}} } \right]^{ - 1} ,\,\,\sigma
h \notin \mathbb{Z},
$$
which completes the proof of Lemma \ref{L5.2}.
\end{proof}

\section{The solution of Problem  3.}

In this section using the results of the previous section we get
explicit formula for coefficients
$\stackrel{\circ}{C}\!\!([\beta];z)$ of the lattice optimal
interpolation formula and we find unknown constant $\lambda$.

Beforehand we give the following result which we use in the proof
of the main theorem.

\begin{lemma}\label{L6.1}
Let $g[\beta]$ be a discrete periodic function, i.e.
$g[\beta]=g(h\beta)=g(h\beta+\gamma)$, $\beta,\gamma\in
\mathbb{Z}$ then the following holds
\begin{equation}
g[\beta]=(g[\beta]\chi_{[0,1]}[\beta])\ast\Phi[\beta],
\label{eq.(6.1)}
\end{equation}
where $\chi_{[0,1]}[\beta]$ is the discrete characteristic
function of the interval $[0,1]$ and $\Phi[\beta]$ is defined by
equality (\ref{eq.(5.4)}).
\end{lemma}

\begin{proof}
  Indeed, using well-known formula  (see. \cite{Sob74})
$$\sum\limits_{k} {\chi
_{[0,1]} \left( {[\beta ] + k} \right) = 1}
$$
and periodicity of the function $g[\beta ]$, taking into account
equality (\ref{eq.(5.4)}), we have
\begin{eqnarray*}
g[\beta ] &=& g[\beta ]\sum\limits_{k} {\chi _{[0,1]} \left(
{[\beta ] + k} \right)}  = \sum\limits_{k} {g\left( {[\beta ] + k}
\right)\chi _{[0,1]} \left( {[\beta ] + k} \right)} \\
 &=& \sum\limits_{k} {g([\beta ] - k)\chi _{[0,1]}
([\beta ] - k) = \sum\limits_{k} {g[\beta  - h^{
- 1} k]\chi _{[0,1]} [\beta  - h^{ - 1} k]  } }\\
&=& \sum\limits_{k} \sum\limits_{\gamma}g[\gamma ]\chi _{[0,1]}
[\gamma ]\delta [\beta  - \gamma  - h^{ - 1} k] \\
&=&\sum\limits_{\gamma}  {g[\gamma ]\chi _{[0,1]} [\gamma
]\sum\limits_{k} {\delta [\beta  - \gamma - h^{ - 1} k]  } }   =
\sum\limits_{\gamma}
{g[\gamma ]\chi _{[0,1]} [\gamma ]\Phi [\beta  - \gamma ]  }\\
&=& \left( {\chi _{[0,1]} [\beta ]g[\beta ]} \right) * \Phi [\beta
].
\end{eqnarray*}
Lemma \ref{L6.1} is proved.
\end{proof}

The main result of the present paper is the following theorem.

\begin{theorem}\label{THM6.2} In the Sobolev space $\widetilde{L_2^{(m)} }(0,1)$
there exists the unique lattice optimal interpolation formula of
the form (\ref{eq.(1.1)}) with the error functional
(\ref{eq.(1.3)}) coefficients of which have the form
\begin{equation}
\stackrel{\circ}{C}\!\!([\beta];z)=h\left(1+\sum\limits_{k\atop
kh\not\in \mathbb{Z}}\frac{\exp(2\pi ik(h\beta-z))}{k^{2m}} \cdot
L(k)\right),\label{eq.(6.2)}
\end{equation}
where $L(k)=\left(\sum\limits_{\scriptstyle
\gamma}\frac{h^{2m}}{|\gamma-hk|^{2m}}\right)^{-1}$.
\end{theorem}

\begin{proof}
Applying the operator $hD_h^{(m)}[\beta ]*$ to both sides of
equation (\ref{eq.(4.3)}) we obtain
 $$
hD_h^{(m)} [\beta ] * \left( {B_{2m} [\beta ] * \left( {\mathop
C\limits^ \circ  ([\beta ],z)\chi _{[0,1]} [\beta ]} \right) +
\lambda } \right) = hD_h^{(m)} [\beta ] * B_{2m} \left( {z -
[\beta ]} \right),\ [\beta ] \in (0,1].
$$
Hence taking into account formulas (\ref{eq.(5.3)}),
(\ref{eq.(6.1)}), (\ref{eq.(5.19)}) we have
$$
\mathop C\limits^ \circ  \left( {[\beta ],z} \right) -
h\sum\limits_{\beta  = 1}^N {\mathop C\limits^ \circ  \left(
{[\beta ],z} \right)}  = hD^{(m)} [\beta ] * B_{2m} (z - [\beta
]),\,\,\,\,[\beta ] \in (0,1].
$$ By virtue of (\ref{eq.(4.4)}) we find
$$
\mathop C\limits^ \circ  \left( {[\beta ],z} \right) = h +
hD^{(m)} [\beta ] * B_{2m} (z - [\beta
]),\,\,\,\,\,\,\,\,\,\,[\beta ] \in (0,1].
$$
Hence using the Bernoulli polynomial
$$
B_{2m} \left( {z - [\beta ]} \right) = \sum\limits_{k \ne 0}
{{{\exp ( - 2\pi ik(z - h\beta ))} \over {(2\pi ik)^{2m} }}}
$$
and equalities  (\ref{eq.(5.18)}), (\ref{eq.(5.19)}) we get
\begin{eqnarray*}
\mathop C\limits^ \circ  \left( {[\beta ],z} \right) &=& h +
hD_h^{(m)} [\beta ] * \sum\limits_{k \ne 0} {{{\exp ( - 2\pi ik(z
-
h\beta ))} \over {(2\pi ik)^{2m} }}}\\
&=& h + h\sum\limits_{k \ne 0} {{{\exp ( - 2\pi ikz)} \over {(2\pi
ik)^{2m} }}} D_h^{(m)} [\beta ] * \exp (2\pi ikh\beta
)\\
&=& h + ( - 1)^m h\sum\limits_{\scriptstyle k \atop
 \scriptstyle kh \notin \mathbb{Z} }
 {{{\exp ( - 2\pi ikz)} \over {(2\pi k)^{2m} }}}
 \cdot ( - 1)^m (2\pi )^{2m} \exp (2\pi ikh\beta ) \left[ {\sum\limits_\gamma  {{{h^{2m} } \over {\left(
{\gamma  - kh} \right)^{2m} }}} } \right]^{ - 1}  \\
&=& h + h\sum\limits_{\scriptstyle k \atop \scriptstyle kh \notin
\mathbb{Z}}{{{(2\pi )^{2m} \exp (2\pi ik(h\beta  - z))} \over
{(2\pi k)^{2m} }}}  \cdot \left[ {\sum\limits_\gamma  {{{h^{2m} }
\over {\left( {\gamma  - kh} \right)^{2m} }}} } \right]^{ - 1}  \\
&=& h + h\sum\limits_{\scriptstyle k \atop \scriptstyle kh \notin
\mathbb{Z}}{{{\exp (2\pi ik(h\beta  - z))} \over {k^{2m} }}} \cdot
\left[ {\sum\limits_\gamma  {{{h^{2m} } \over {\left( {\gamma  -
kh} \right)^{2m} }}} } \right]^{ - 1}.
\end{eqnarray*}
Hence, setting
\begin{equation}
L(k) = \left[ {\sum\limits_\gamma {{{h^{2m} } \over {\left(
{\gamma  - kh} \right)^{2m} }}} } \right]^{ - 1},\label{eq.Lk}
\end{equation}  we get
(\ref{eq.(6.2)}). Theorem \ref{THM6.2} is proved.
\end{proof}

Now using Theorem \ref{THM6.2} we find $\lambda$. First, putting
the expression (\ref{eq.(6.2)}) of the coefficients  $\mathop
C\limits^ \circ \left( {[\beta ],z} \right)$ to the left side of
equality (\ref{eq.(4.1)}) we calculate the following sum
\begin{eqnarray*}
Q(z) &=& \sum\limits_{\gamma  = 1}^N {\mathop C\limits^ \circ
([\gamma ],z)B_{2m} [\beta  - \gamma ] = }\\
 &=& \sum\limits_{\gamma  = 1}^N {h\left( {1 + \sum\limits_{\scriptstyle k \atop
\scriptstyle kh \notin Z} ^{} {{{\exp (2\pi ik(h\gamma  - z)L(k)}
\over {k^{2m} }}} } \right)\left( {\sum\limits_{\alpha  \ne 0}^{}
{{{\exp ( - 2\pi i\alpha (h\beta  - h\gamma )} \over {(2\pi
i\alpha )^{2m} }}} } \right)  }\\
&=& \sum\limits_{\alpha  \ne 0}^{} {{{\exp ( - 2\pi i\alpha h\beta
)} \over {(2\pi i\alpha )^{2m} }}\sum\limits_{\gamma  = 1}^N
{h\exp (2\pi i\alpha h\gamma )}   }\\
&& + \sum\limits_{\scriptstyle k \atop \scriptstyle kh \notin Z}
^{} {{{\exp ( - 2\pi ikz)\, \cdot L(k)} \over {k^{2m} }}}
\,\sum\limits_{\alpha  \ne 0}^{} {{{\exp ( - 2\pi i\alpha h\beta
)} \over {(2\pi i\alpha )^{2m} }}} \cdot \sum\limits_{\gamma  =
1}^N {h\exp (2\pi ikh\gamma  + 2\pi i\alpha h\gamma ).}
\end{eqnarray*}
It is known that
\begin{eqnarray}
\sum\limits_{\gamma  = 1}^N {h\exp (2\pi i\alpha h\gamma )} & =&
{{h\exp (2\pi i\alpha h) \cdot (1 - \exp (2\pi i\alpha ))} \over
{1 - \exp (2\pi i\alpha h)}} \nonumber\\
 &=& \left\{ \begin{array}{l}
  0,\,\,\alpha h \notin \mathbb{Z}, \cr
  1,\,\,\alpha h \in \mathbb{Z},  \cr
 \end{array}
 \right.\label{eq.(6.3)}
\end{eqnarray}
\begin{eqnarray}
\sum\limits_{\gamma  = 1}^N h\exp (2\pi i(k + \alpha )h\gamma )
&=& {{h\exp (2\pi i(k + \alpha )h) \cdot \left( {1 - \exp (2\pi
i(k + \alpha ))} \right)} \over {1 - \exp (2\pi i(k + \alpha )h)}}  \nonumber \\
 &=& \left\{ \begin{array}{l}
  0,\,\,(k + \alpha )h \notin \mathbb{Z}, \hfill \cr
  1,\,\,(k + \alpha )h \in \mathbb{Z}.  \cr
  \end{array}  \right.
  \label{eq.(6.4)}
\end{eqnarray}
By virtue of  equalities (\ref{eq.(6.3)}), (\ref{eq.(6.4)}) the
expression $Q(z)$ takes the form
    $$
Q(z) = \sum\limits_{\scriptstyle \gamma  \ne 0 \hfill \atop
\scriptstyle \gamma h \in \mathbb{Z}}{{1 \over {(2\pi i\gamma
)^{2m} }} + } \sum\limits_{\scriptstyle k \hfill \atop
 \scriptstyle kh \notin \mathbb{Z}}{{{L(k)} \over {(2\pi ik)^{2m} }} }
\sum\limits_{\alpha  \ne 0}^{} {{{\exp ( - 2\pi i(\alpha h\beta  +
kz)} \over {\alpha ^{2m} }}},
$$
when  $(k + \alpha )h \in \mathbb{Z}.$\\
From here setting
  $$
(k + \alpha )h = t,\,\,\,t \in \mathbb{Z},\,\,\,k + \alpha  = th^{
- 1} ,\,\,\,\,\alpha  = th^{ - 1}  - k,
$$
we have
\begin{eqnarray*}
Q(z) &=& \sum\limits_{\scriptstyle \gamma  \ne 0 \hfill \atop
  \scriptstyle \gamma h \in \mathbb{Z} } {{1 \over {(2\pi i\gamma )^{2m} }}  }
+ \sum\limits_{\scriptstyle k \hfill \atop \scriptstyle kh \notin
\mathbb{Z} } {{{L(k)} \over {(2\pi ik)^{2m} }}\sum\limits_t^{}
{{{\exp ( - 2\pi i((th^{ - 1}  - k)h\beta  + kz)} \over {(th^{ -
1}  - k)^{2m} }}}  }\\
&=& \sum\limits_{\scriptstyle \gamma  \ne 0 \hfill \atop
\scriptstyle \gamma h \in \mathbb{Z}} {{1 \over {(2\pi i\gamma
)^{2m} }} + \sum\limits_{\scriptstyle k \hfill \atop \scriptstyle
kh \notin \mathbb{Z }}{{{L(k)} \over {(2\pi ik)^{2m}
}}\sum\limits_t^{} {{{\exp ( - 2\pi i((t\beta  - kh\beta  + kz)}
\over {(th^{ - 1}  - k)^{2m} }}  } } }\\
&=& \sum\limits_{\scriptstyle \gamma  \ne 0 \hfill \atop
\scriptstyle \gamma h \in \mathbb{Z}}  {{1 \over {(2\pi i\gamma
)^{2m} }} + \sum\limits_{\scriptstyle k \hfill \atop \scriptstyle
kh \notin \mathbb{Z}} {{{\exp ( - 2\pi ik(z - h\beta ) \cdot L(k)}
\over {(2\pi ik)^{2m} }}\sum\limits_t^{} {{{\exp ( - 2\pi it\beta
)} \over {(th^{ - 1}  - k)^{2m} }}.} } }
\end{eqnarray*}
Hence keeping in mind that $\exp(-2\pi it\beta ) = 1$ and $L(k) =
\sum\limits_t^{} {{1 \over {(th^{ - 1}  - k)^{2m} }}} $ we obtain
\begin{equation}
Q(z) = \sum\limits_{\scriptstyle \gamma  \ne 0 \hfill \atop
\scriptstyle \gamma h \in \mathbb{Z}}  {{1 \over {(2\pi i\gamma
)^{2m} }} + \sum\limits_{\scriptstyle k \hfill \atop \scriptstyle
kh \notin \mathbb{Z}}{{{\exp ( - 2\pi ik(z - h\beta )} \over
{(2\pi ik)^{2m} }}}.}\label{eq.(6.5)}
\end{equation}
Adding and subtracting to the right hand side of equality
(\ref{eq.(6.5)}) the following series
 $$
\sum\limits_{\scriptstyle k \ne 0 \hfill \atop \scriptstyle kh \in
\mathbb{Z }}{{{\exp ( - 2\pi ik(z - h\beta ))} \over {(2\pi
ik)^{2m} }}},
$$
we have
\begin{eqnarray}
 Q(z) &=& \sum\limits_{\scriptstyle \gamma  \ne 0 \hfill \atop
\scriptstyle \gamma h \in \mathbb{Z}}{{1 \over {(2\pi i\gamma
)^{2m} }} - \sum\limits_{\scriptstyle k \ne 0 \hfill \atop
\scriptstyle kh \in \mathbb{Z} }{{{\exp (2\pi ik(h\beta  - z)}
\over {(2\pi ik)^{2m} }} + B_{2m} (z - h\beta )  } }\nonumber\\
& =& \sum\limits_{\scriptstyle k \ne 0 \hfill \atop \scriptstyle
kh \in \mathbb{Z}}{{{1 - \exp ( - 2\pi ikz)} \over {(2\pi ik)^{2m}
}} + B_{2m} (z - h\beta )}.\label{eq.(6.6)}
\end{eqnarray}
Putting the expression (\ref{eq.(6.6)}) of $Q(z)$ to the left hand
side of equation (\ref{eq.(4.1)}) for $\lambda $ we obtain the
following expression
 \begin{equation}
\lambda  = \sum\limits_{\scriptstyle k \ne 0 \hfill \atop
\scriptstyle kh \in \mathbb{Z}}{{{\exp ( - 2\pi ikz) - 1} \over
{(2\pi ik)^{2m} }}.}\label{eq.(6.7)}
 \end{equation}
Hence clear that when $z = h\beta$
$$
\lambda  = 0.
$$

Now we show that the expression (\ref{eq.(6.2)}) of the
coefficients $\mathop C\limits^ \circ \left( {[\beta ],z} \right)$
satisfy equality (\ref{eq.(4.2)}). So, we have
\begin{eqnarray*}
\sum\limits_{\beta  = 1}^N {\mathop C\limits^ \circ  \left(
{[\beta ],z} \right)}  &=& \sum\limits_{\beta  = 1}^N {\left( {h +
h\sum\limits_{k \atop kh \notin \mathbb{Z}} {{{\exp (2\pi
ik(h\beta - z))} \over
{k^{2m} }}}  \cdot L(k)} \right)}\\
& =& \sum\limits_{\beta  = 1}^N h  + h\sum\limits_{k \atop kh
\notin \mathbb{Z}} {{{\exp ( - 2\pi ikz)L(k)} \over {k^{2m} }}}
\sum\limits_{\beta  =
1}^N {\exp (2\pi ikh\beta )}  \\
&=& Nh + h\sum\limits_{k \atop kh \notin \mathbb{Z}} {{{\exp ( -
2\pi ikz)L(k)} \over {k^{2m} }}}  \cdot {{\exp (2\pi ikh)(1 -
e^{2\pi ik} )}
\over {1 - e^{2\pi ikh} }}\\
 &=& 1.
\end{eqnarray*}

Thus Problem 3 and respectively Problem 2 are solved.

\section{The norm of the error functional\\ of lattice optimal
interpolation formulas}

In this section using the results of previous sections we
calculate the norm of the error functional $\ell$ of the lattice
optimal interpolation formula.

The following holds

\begin{theorem}\label{THM7.1}
Square of the norm of the error functional (\ref{eq.(1.3)}) of the
lattice optimal interpolation formula of the form (\ref{eq.(1.1)})
in the space $\widetilde{L_2^{(m)}}(0,1)$ have the following form
\begin{eqnarray}
 \left\| {\mathop \ell \limits^ \circ }
|\widetilde{L_2^{(m)}}(0,1)\right\|^2 & =& ( - 1)^m \Bigg[
{\sum\limits_{k \ne 0}^{} {{1 \over {(2\pi ik)^{2m} }} -
\sum\limits_{\scriptstyle k \ne 0 \hfill \atop \scriptstyle kh \in
\mathbb{Z}}{{{2\exp ( - 2\pi ikz) - 1}
\over {(2\pi ik)^{2m} }} } } } \nonumber \\
&& - \sum\limits_{\scriptstyle k \hfill \atop \scriptstyle kh
\notin \mathbb{Z}}{{{L(k)} \over {(2\pi ik)^{2m}
}}\sum\limits_t^{} { {{{\exp (2\pi izth^{ - 1} )} \over {(th^{ -
1}  - k)^{2m} }}} } }\Bigg],\label{eq.(7.1)}
\end{eqnarray}
where $L(k)$ is defined by equality (\ref{eq.Lk}).
\end{theorem}

\begin{proof}
From equality (\ref{eq.(2.7)}) when $x_k  = hk$ for $\left\|
{\mathop \ell \limits^ \circ  } \right\|^2$ we have
\begin{eqnarray*}
\left\| {\mathop \ell \limits^ \circ  } \right\|^2  &=& ( - 1)^m
\left[ {B_{2m} (0) - 2\sum\limits_{\gamma  = 1}^N {\mathop
C\limits^ \circ  \left( {[\gamma ],z} \right)B_{2m} (z - h\gamma )
+ } } \right.\\
&&\left. { + \sum\limits_{\gamma  = 1}^N {\mathop C\limits^ \circ
} ([\gamma ],z)\sum\limits_{\beta  = 1}^N {\mathop C\limits^ \circ
([\beta ],z)B_{2m} (h\beta  - h\gamma )} } \right].
\end{eqnarray*}
Hence taking into account (\ref{eq.(4.1)}) and (\ref{eq.(4.2)}) we
obtain
\begin{eqnarray*}
\left\| {\mathop \ell \limits^ \circ  } \right\|^2  &=& ( - 1)^m
\left[ {B_{2m} (0) - 2\sum\limits_{\gamma  = 1}^N {\mathop
C\limits^ \circ  \left( {[\gamma ],z} \right)B_{2m} (z - h\gamma )
 } } \right. \left. { + \sum\limits_{\beta  = 1}^N {\mathop
C\limits^ \circ  } ([\beta ],z)\left( {B_{2m} (z - h\beta ) -
\lambda } \right)} \right]\\
& =&(-1)^m \left[ {B_{2m} (0) - \sum\limits_{\gamma  = 1}^N
{\mathop C\limits^ \circ  ([\gamma ],z)B_{2m} (z - h\gamma ) -
\lambda } } \right].
\end{eqnarray*}
From here using the expression (\ref{eq.(6.7)}) of $\lambda$ we
obtain
\begin{eqnarray}
\left\| {\mathop \ell \limits^ \circ } \right\|^2 & =& ( - 1)^m
\Bigg[ {B_{2m} (0) - \sum\limits_{\gamma = 1}^N {h\Bigg( {1 +
\sum\limits_{\scriptstyle k \hfill \atop \scriptstyle hk \notin
\mathbb{Z} } {{{\exp (2\pi ik(h\gamma - z)\,L(k)} \over {k^{2m}
}}} } \Bigg) } }\nonumber \\
&& \times  {\sum\limits_{\alpha  \ne 0}^{} {{{\exp ( - 2\pi
i\alpha (z - h\gamma ))} \over {(2\pi i\alpha )^{2m} }} -
\sum\limits_{\scriptstyle k \ne 0 \hfill \atop \scriptstyle kh \in
\mathbb{Z} }{{{\exp ( - 2\pi izk) - 1} \over {(2\pi ik)^{2m} }}} }
} \Bigg] \nonumber\\
&=&(-1)^{2m} \Bigg[ {\sum\limits_{k \ne 0}^{} {{1 \over {(2\pi
ik)^{2m} }} - \sum\limits_{\alpha  \ne 0}^{} {{{\exp ( - 2\pi
i\alpha z)} \over {(2\pi i\alpha )^{2m} }}  } } }
\sum\limits_{\gamma  = 1}^N h\exp (2\pi i\alpha h\gamma )\nonumber \\
&& - \sum\limits_{\scriptstyle k \hfill \atop \scriptstyle kh
\notin \mathbb{Z}}{{{\exp ( - 2\pi ikz)\,L(k)} \over {k^{2m} }} }
 \sum\limits_{\alpha  \ne 0} {{{\exp ( - 2\pi i\alpha z)}
\over {(2\pi i\alpha )^{2m} }}} \sum\limits_{\gamma  = 1}^N {h\exp
(2\pi i(\alpha  + k)h\gamma )  } \nonumber\\
&&-\sum\limits_{\scriptstyle k \ne 0 \hfill \atop \scriptstyle kh
\in Z \hfill} {{{\exp ( - 2\pi ikz) - 1} \over {(2\pi ik)^{2m} }}}
\Bigg].\label{eq.(7.2)}
\end{eqnarray}
From (\ref{eq.(7.2)}) taking into account equalities
(\ref{eq.(6.3)}) and (\ref{eq.(6.4)}) we obtain
\begin{eqnarray*}
\left\| {\mathop \ell \limits^ \circ } \right\|^2 & =& ( - 1)^m
\Bigg[ {\sum\limits_{k \ne 0}^{} {{1 \over {(2\pi ik)^{2m} }} -
\sum\limits_{\scriptstyle k  \ne 0 \atop \scriptstyle kh \in
\mathbb{Z} } {{{\exp ( - 2\pi ik z)}
\over {(2\pi ik )^{2m} }}  } } } \\
&&- \sum\limits_{\scriptstyle k \atop \scriptstyle kh \notin
\mathbb{Z}} {{{L(k)} \over {k^{2m} }}\sum\limits_{\scriptstyle
\alpha  \ne 0  \atop \scriptstyle (k + \alpha )h \in
\mathbb{Z}}{{{\exp ( - 2\pi iz(k + \alpha ))} \over {(2\pi i\alpha
)^{2m} }}-} } {\sum\limits_{\scriptstyle k \ne 0 \hfill \atop
\scriptstyle kh \in \mathbb{Z}}{{{\exp ( - 2\pi ikz) - 1} \over
{(2\pi ik)^{2m} }}} } \Bigg].
\end{eqnarray*}
Whence, setting $(k+\alpha)h=t$, $\alpha=th^{-1}-k$, we get
(\ref{eq.(7.1)}). Theorem \ref{THM7.1} is proved.
\end{proof}

Now we show that the expression (\ref{eq.(7.1)}) is zero at the
nodes $x_{\beta} = h\beta$ of the lattice optimal interpolation
formula.

Suppose $z=h\beta$, then from (\ref{eq.(7.1)}) we obtain
\begin{eqnarray*}
 \left\| {\mathop \ell \limits^ \circ  } \right\|^2  &=& ( - 1)^m
\left[ {\sum\limits_{k \ne 0}^{} {{1 \over {(2\pi ik)^{2m} }} -
\sum\limits_{\scriptstyle k \ne 0 \hfill \atop  \scriptstyle kh
\in \mathbb{Z} } {{{2\exp ( - 2\pi ikh\beta ) - 1} \over {(2\pi
ik)^{2m} }}  } } } \right. - \left. {\sum\limits_{\scriptstyle k
\hfill \atop \scriptstyle kh \notin \mathbb{Z} }{{{L(k)} \over
{(2\pi ik)^{2m} }}\sum\limits_t^{} {{{\exp (2\pi it\beta )} \over
{(th^{ - 1}  - k)^{2m} }}} } } \right] \\
&=& ( - 1)^m \left[ {\sum\limits_{k \ne 0}^{} {{1 \over {(2\pi
ik)^{2m} }} - \sum\limits_{\scriptstyle k \ne 0 \hfill \atop
 \scriptstyle kh \in \mathbb{Z} }{{1 \over {(2\pi ik)^{2m} }} - \left. {\sum\limits_{\scriptstyle k \hfill \atop
\scriptstyle kh \notin \mathbb{Z} }{{{L(k)} \over {(2\pi ik)^{2m}
\cdot L(k)}}} } \right]} } } \right. \\
&=& ( - 1)^m \left[ {\sum\limits_{k \ne 0}^{} {{1 \over {(2\pi
ik)^{2m} }} - \sum\limits_{\scriptstyle k \ne 0 \hfill \atop
\scriptstyle kh \in \mathbb{Z} } {{1 \over {(2\pi ik)^{2m} }} -
\sum\limits_{\scriptstyle k \hfill \atop \scriptstyle kh \notin
\mathbb{Z} }{{1 \over {(2\pi ik)^{2m} }}} } } } \right]\\
& =& 0.
\end{eqnarray*}

This means that the condition of interpolation is fulfilled and
this confirms our theoretical results.

\section{Connection between lattice optimal interpolation
formula\\
and optimal quadrature formula in the space
$\widetilde{L_2^{(m)}}(0,1)$}

In previous sections we constructed the lattice optimal
interpolation formula which has the following form
\begin{equation}
\label{eq.(8.1)} \varphi(x)\cong
\sum\limits_{\beta=1}^N\stackrel{\circ}{C}\!([\beta],x)\varphi(h\beta),
\end{equation}
where $\stackrel{\circ}{C}\!([\beta],x)$ are defined by expression
(\ref{eq.(6.2)}), $h=\frac{1}{N},$ $N=2,3,...$.

Integrating equality (\ref{eq.(8.1)}) from 0 to 1, we get
\begin{eqnarray}
\int_0^1\varphi(x)dx&\cong&
\int_0^1\sum\limits_{\beta=1}^N\stackrel{\circ}{C}([\beta],x)\varphi(h\beta)dx
\nonumber \\
&\cong&
\sum\limits_{\beta=1}^N\left(\int_0^1h\varphi(h\beta)dx+h\int_0^1\sum\limits_{k\atop
kh\notin \mathbb{Z}}\displaystyle \frac{\exp(2\pi i
k(h\beta-x))\cdot L(k)}{k^{2m}}\right).\label{eq.(8.2)}
\end{eqnarray}
The second integral in the right hand side of (\ref{eq.(8.2)}) is
equal to zero, since
$\int\limits_0^1\exp(2\pi i k(h\beta-x)dx=0$.\\
Then from (\ref{eq.(8.2)}) we get the following  well-known
quadrature formula
$$
\int_0^1\varphi(x)dx\cong h\sum_{\beta=1}^N\varphi(h\beta).
$$
This is the rectangular formula and optimality of this quadrature
formula in the space $\widetilde{L_2^{(m)}}(0,1)$ is known (see,
for example \cite{SMNikBook88,Sob74}).

Thus by integrating the optimal interpolation formula in the space
  $\widetilde{L_2^{(m)}}(0,1)$ we have obtained the optimal quadrature
  formula in the same space.

\section{Numerical results}

First of all here we give explicit formulas for coefficients of
lattice optimal interpolation formula (\ref{eq.(1.1)}) which are
very useful in practice. From Theorem \ref{THM6.2} we get the
following

\begin{corollary}
\label{Cor9.1} In the Sobolev space $\widetilde{L_2^{(m)} }(0,1)$
there exists the unique lattice optimal interpolation formula of
the form (\ref{eq.(1.1)}) with the error functional
(\ref{eq.(1.3)}) coefficients of which have the form
\begin{equation}
\stackrel{\circ}{C}\!\!([\beta];z)=h\Bigg[1+\sum\limits_{j=1}^{\frac{N-1}{2}}\sum\limits_t
\frac{2\cos(2\pi(Nt+j)(h\beta-z))}{(Nt+j)^{2m}}\left(\sum\limits_{\gamma}\frac{1}{(N\gamma+j)^{2m}}\right)^{-1}\Bigg],\
\beta=\overline{1,N}\label{eq.(9.1)}
\end{equation}
when $N$ is odd number and $N\geq m$ and
\begin{eqnarray}
\stackrel{\circ}{C}\!\!([\beta];z)&=&h\Bigg[1+\sum\limits_{j=1}^{\frac{N}{2}-1}\sum\limits_t
\frac{2\cos(2\pi(Nt+j)(h\beta-z))}{(Nt+j)^{2m}}\left(\sum\limits_{\gamma}\frac{1}{(N\gamma+j)^{2m}}\right)^{-1}\nonumber\\
&&
+\sum\limits_t\frac{\cos(2\pi(Nt+N/2)(h\beta-z))}{(Nt+N/2)^{2m}}\left(\sum\limits_{\gamma}\frac{1}{(N\gamma+N/2)^{2m}}\right)^{-1}
\Bigg],\ \beta=\overline{1,N}\label{eq.(9.2)}
\end{eqnarray} when
$N$ is even number and $N\geq m$.
\end{corollary}

\begin{proof} Now we simplify the expression (\ref{eq.(6.2)}) of
the optimal coefficients $\stackrel{\circ}{C}\!\!([\beta];z)$.

We denote
\begin{equation}\label{eq.(9.3)}
S=\sum\limits_{k\atop kh\not\in \mathbb{Z}}\frac{\exp(2\pi
ik(h\beta-z))}{k^{2m}} \cdot \left(\sum\limits_{\scriptstyle
\gamma}\frac{1}{|N\gamma-k|^{2m}}\right)^{-1}
\end{equation}

Consider two cases.

\emph{The case 1.} Let $N$ be odd number, i.e. $N=2M+1$ then
$k\neq Nt, \ t\in \mathbb{Z}$. This means
$$
k=\left\{
\begin{array}{l}
Nt+1,\ Nt+2, ...,Nt+M,\\
Nt-1,\ Nt-2, ...,Nt-M.\\
\end{array}
\right.
$$
Therefore from (\ref{eq.(9.3)}) we have
\begin{eqnarray*}
S&=&\sum\limits_t\sum\limits_{j=1}^M\Bigg[\frac{\exp(2\pi
i(Nt+j)(h\beta-z))}{(Nt+j)^{2m}}\left(\sum\limits_{\gamma}\frac{1}{(N\gamma-(Nt+j))^{2m}}\right)^{-1}\\
&& + \frac{\exp(2\pi
i(Nt-j)(h\beta-z))}{(Nt-j)^{2m}}\left(\sum\limits_{\gamma}\frac{1}{(N\gamma-(Nt-j))^{2m}}\right)^{-1}\Bigg].
\end{eqnarray*}
Hence, taking into account that $\gamma$ and $t$ take all integer
value, replacing $\gamma$ by $-\gamma$ and $t$ by $-t$ in the
second term of the square brackets we get
\begin{eqnarray*}
S&=&\sum\limits_t\sum\limits_{j=1}^M\Bigg[\frac{\exp(2\pi
i(Nt+j)(h\beta-z))}{(Nt+j)^{2m}} + \frac{\exp(-2\pi
i(Nt+j)(h\beta-z))}{(Nt+j)^{2m}}\Bigg]\left(\sum\limits_{\gamma}\frac{1}{(N\gamma+j)^{2m}}\right)^{-1}\\
&=&\sum\limits_{j=1}^{M}\sum\limits_t
\frac{2\cos(2\pi(Nt+j)(h\beta-z))}{(Nt+j)^{2m}}\left(\sum\limits_{\gamma}\frac{1}{(N\gamma+j)^{2m}}\right)^{-1}.
\end{eqnarray*}
Taking into account that $M=(N-1)/2$ and putting the last equality
for $S$ into (\ref{eq.(6.2)}) we get (\ref{eq.(9.1)}).

\emph{The case 2.} Let $N$ be even number, i.e. $N=2M$ then $k\neq
Nt, \ t\in \mathbb{Z}$. This means
$$
k=\left\{
\begin{array}{l}
Nt+1,\ Nt+2, ...,Nt+M-1,Nt+M,\\
Nt-1,\ Nt-2, ...,Nt-M+1.\\
\end{array}
\right.
$$
Therefore from (\ref{eq.(9.3)}) we have
\begin{eqnarray*}
S&=&\sum\limits_t\sum\limits_{j=1}^{M-1}\Bigg[\frac{\exp(2\pi
i(Nt+j)(h\beta-z))}{(Nt+j)^{2m}}\left(\sum\limits_{\gamma}\frac{1}{(N\gamma-(Nt+j))^{2m}}\right)^{-1}\\
&& + \frac{\exp(2\pi
i(Nt-j)(h\beta-z))}{(Nt-j)^{2m}}\left(\sum\limits_{\gamma}\frac{1}{(N\gamma-(Nt-j))^{2m}}\right)^{-1}\Bigg]\\
&&+\sum\limits_t\frac12\Bigg[\frac{\exp(2\pi
i(Nt+M)(h\beta-z))}{(Nt+M)^{2m}}+\frac{\exp(2\pi
i(Nt-M)(h\beta-z))}{(Nt-M)^{2m}}\Bigg]\left(\sum\limits_{\gamma}\frac{1}{(N\gamma+M)^{2m}}\right)^{-1}.
\end{eqnarray*}
Hence as in the case 1, taking into account that $\gamma$ and $t$
take all integer value, replacing $\gamma$ by $-\gamma$ and $t$ by
$-t$ in the second terms of the square brackets we get
\begin{eqnarray*}
S&=&\sum\limits_t\sum\limits_{j=1}^{M-1}\Bigg[\frac{\exp(2\pi
i(Nt+j)(h\beta-z))}{(Nt+j)^{2m}} + \frac{\exp(-2\pi
i(Nt+j)(h\beta-z))}{(Nt+j)^{2m}}\Bigg]\left(\sum\limits_{\gamma}\frac{1}{(N\gamma+j)^{2m}}\right)^{-1}\\
&&+\sum\limits_t\frac12\Bigg[\frac{\exp(2\pi
i(Nt+M)(h\beta-z))}{(Nt+M)^{2m}}+\frac{\exp(-2\pi
i(Nt+M)(h\beta-z))}{(Nt+M)^{2m}}\Bigg]\left(\sum\limits_{\gamma}\frac{1}{(N\gamma+M)^{2m}}\right)^{-1}\\
&=&\sum\limits_{j=1}^{M-1}\sum\limits_t
\frac{2\cos(2\pi(Nt+j)(h\beta-z))}{(Nt+j)^{2m}}\left(\sum\limits_{\gamma}\frac{1}{(N\gamma+j)^{2m}}\right)^{-1}\\
&& +\sum\limits_t
\frac{\cos(2\pi(Nt+M)(h\beta-z))}{(Nt+M)^{2m}}\left(\sum\limits_{\gamma}\frac{1}{(N\gamma+M)^{2m}}\right)^{-1}.
\end{eqnarray*}
Taking into account that $M=N/2$ and putting the last equality for
$S$ into (\ref{eq.(6.2)}) we get (\ref{eq.(9.2)}).\\ Corollary
\ref{Cor9.1} is proved.
\end{proof}

\newpage

Here in numerical examples we consider interpolation of two 1-periodic functions:\\
1) $\varphi_1(x)=\sin(2\pi x)$,\\ 2) the Bernoulli polynomial of
degree $10$, i.e.\\
$$\varphi_2(x)\equiv
B_{10}(x)=x^{10}-5x^9+\frac{15}{2}x^8-7x^6+5x^4-\frac32
x^2+\frac{5}{66}.
$$

We interpolate these two functions with optimal interpolation
formula
\begin{equation}
\label{eq.(9.4)}
P_{\varphi}(x)=\sum\limits_{\beta=1}^5\stackrel{\circ}{C}\!\!([\beta];x)\varphi(h\beta)
\end{equation}
for $m=1,2,3,4$. For simplicity we have taken $N=5$. Since in our
case $N=5$, $h=1/5$ as optimal coefficients we use formula
(\ref{eq.(9.1)}). One can use formula (\ref{eq.(9.2)}) for even
number $N$ of the nodes.

Below in each case graphs of the optimal coefficients and graphs
of absolute errors between optimal interpolation formula
(\ref{eq.(9.4)}) and functions $\sin(2\pi x)$, $B_{10}(x)$
 are respectively given.

\emph{The case $m=1$.}
\begin{figure}[h]
\includegraphics[width=60pt,angle=-90]{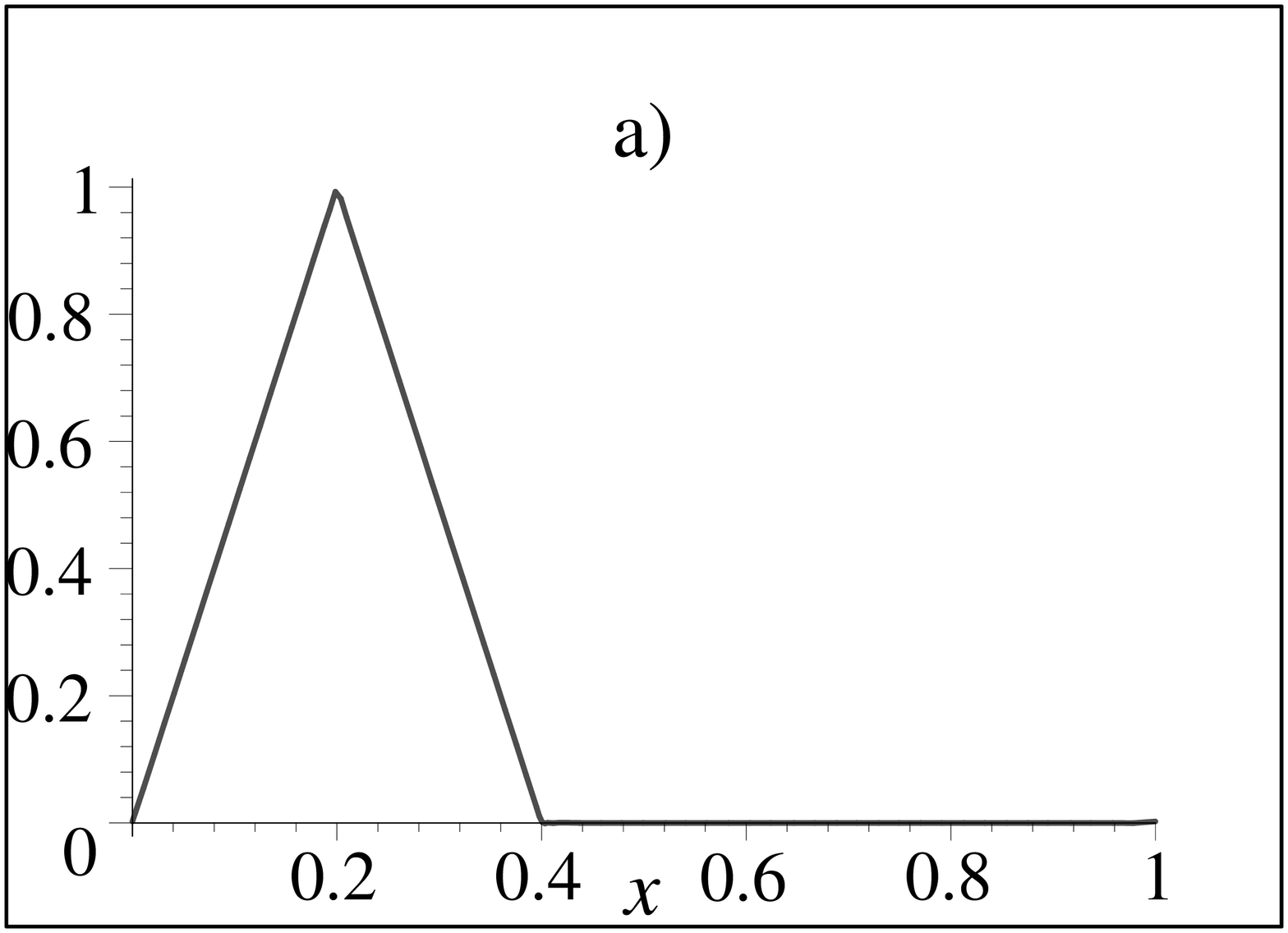}
\includegraphics[width=60pt,angle=-90]{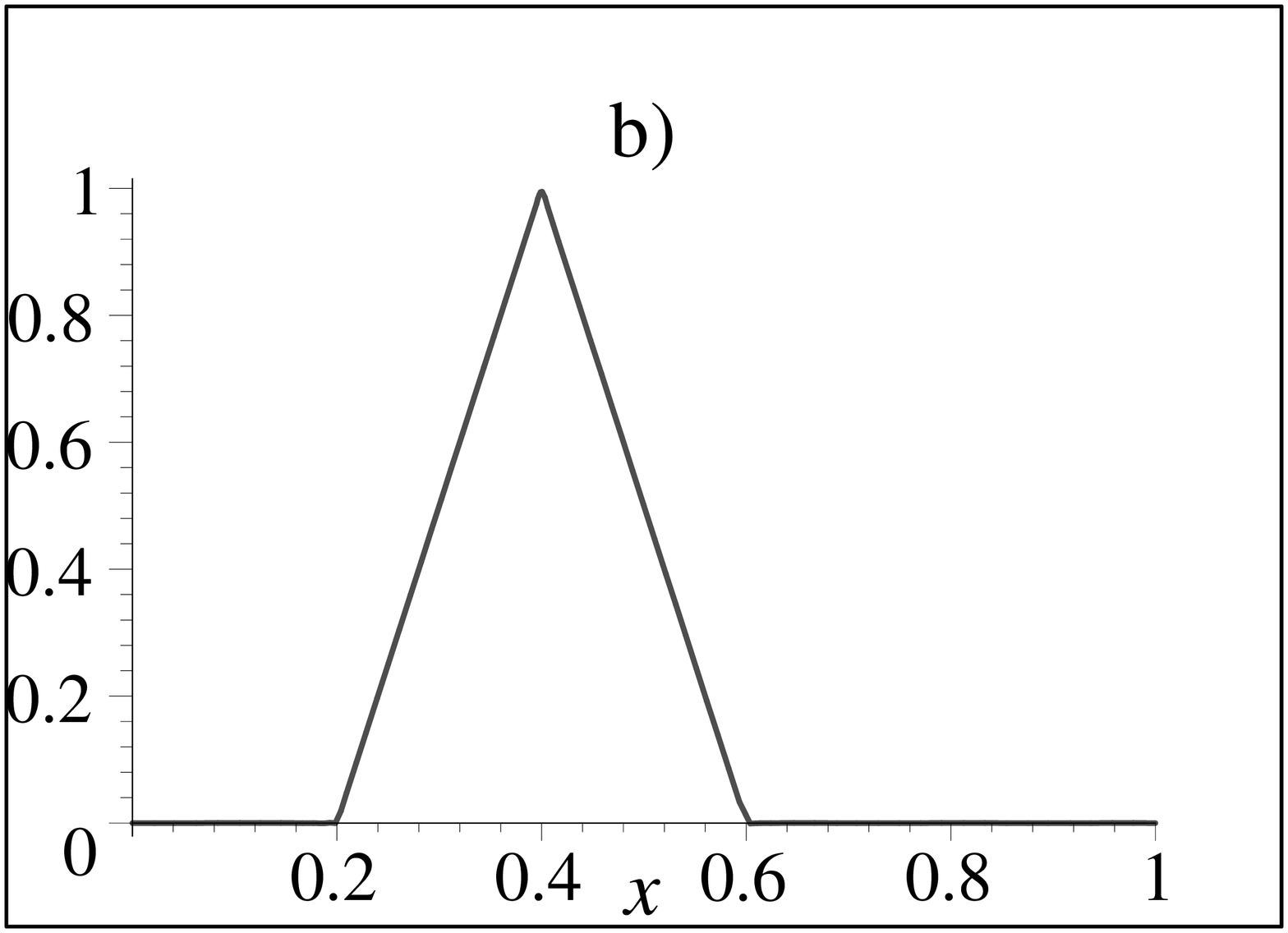}
\includegraphics[width=60pt,angle=-90]{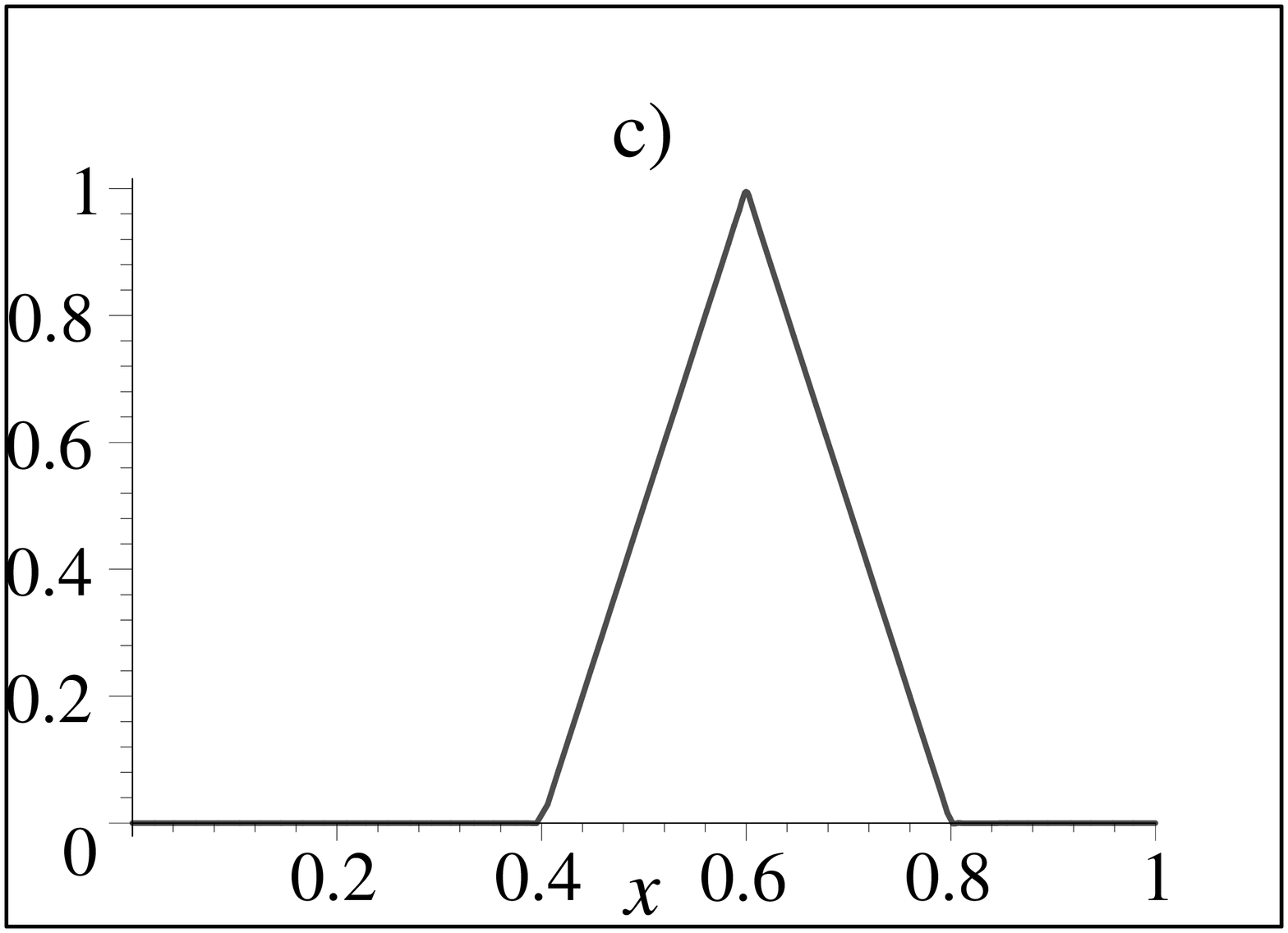}
\includegraphics[width=60pt,angle=-90]{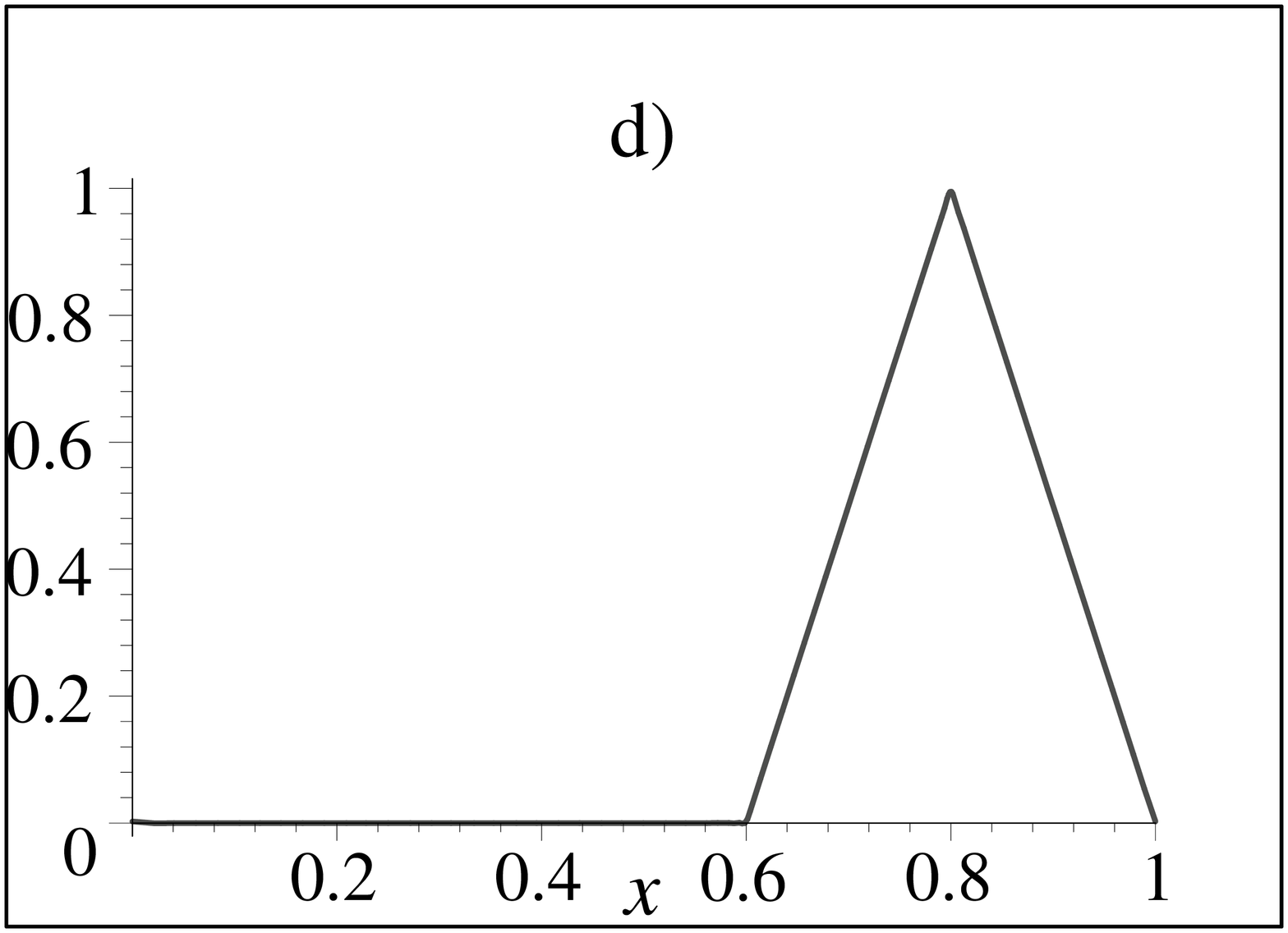}
\includegraphics[width=60pt,angle=-90]{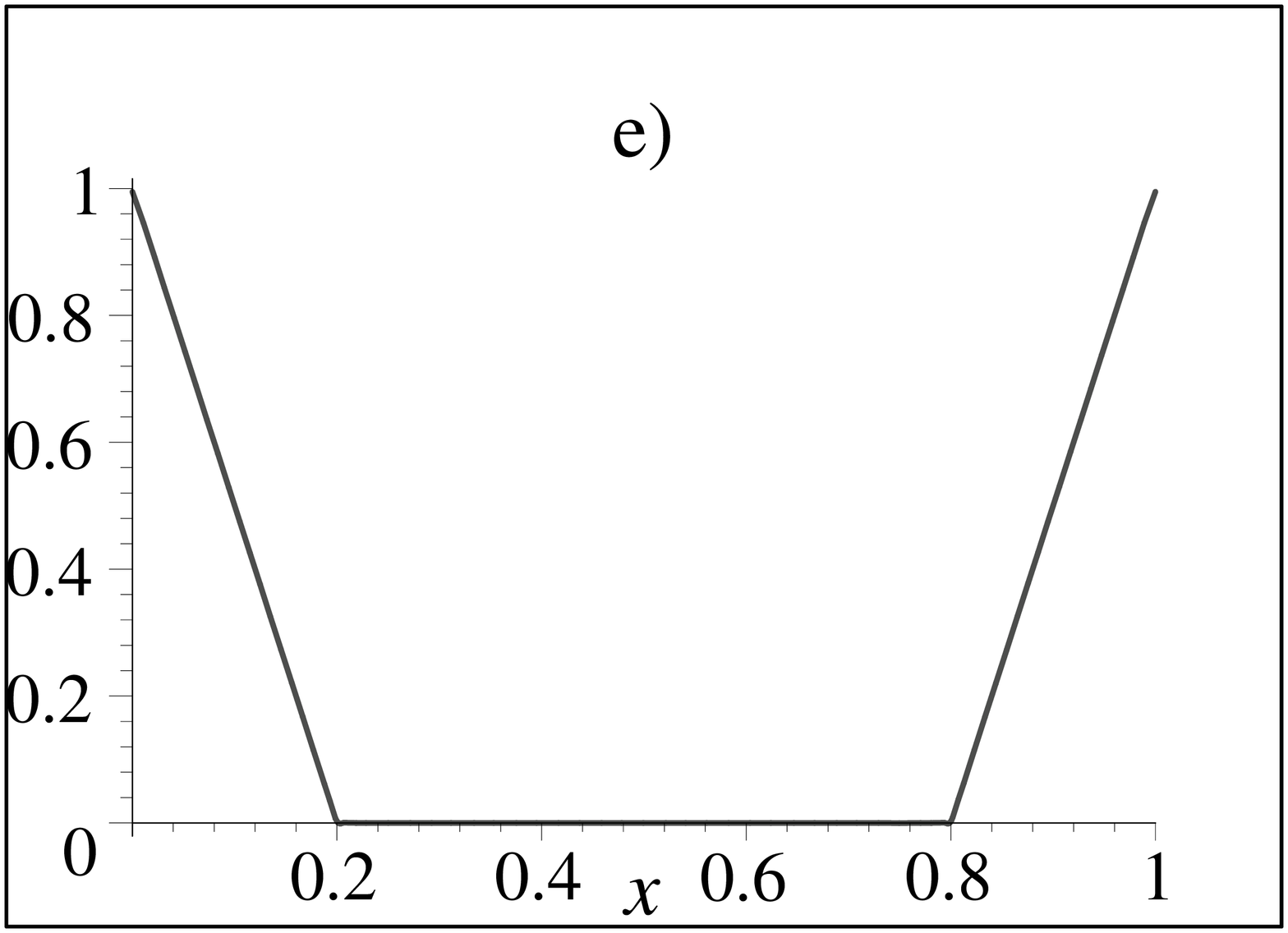}\\
\caption{ Graphs of optimal coefficients for $m=1$, $N=5$:
 a) $C([1],x)$, b) $C([2],x)$, c) $C([3],x)$, d) $C([4],x)$, e)
$C([5],x)$. } \label{Fig1}
\end{figure}
\begin{figure}[h]
\includegraphics[width=100pt,angle=-90]{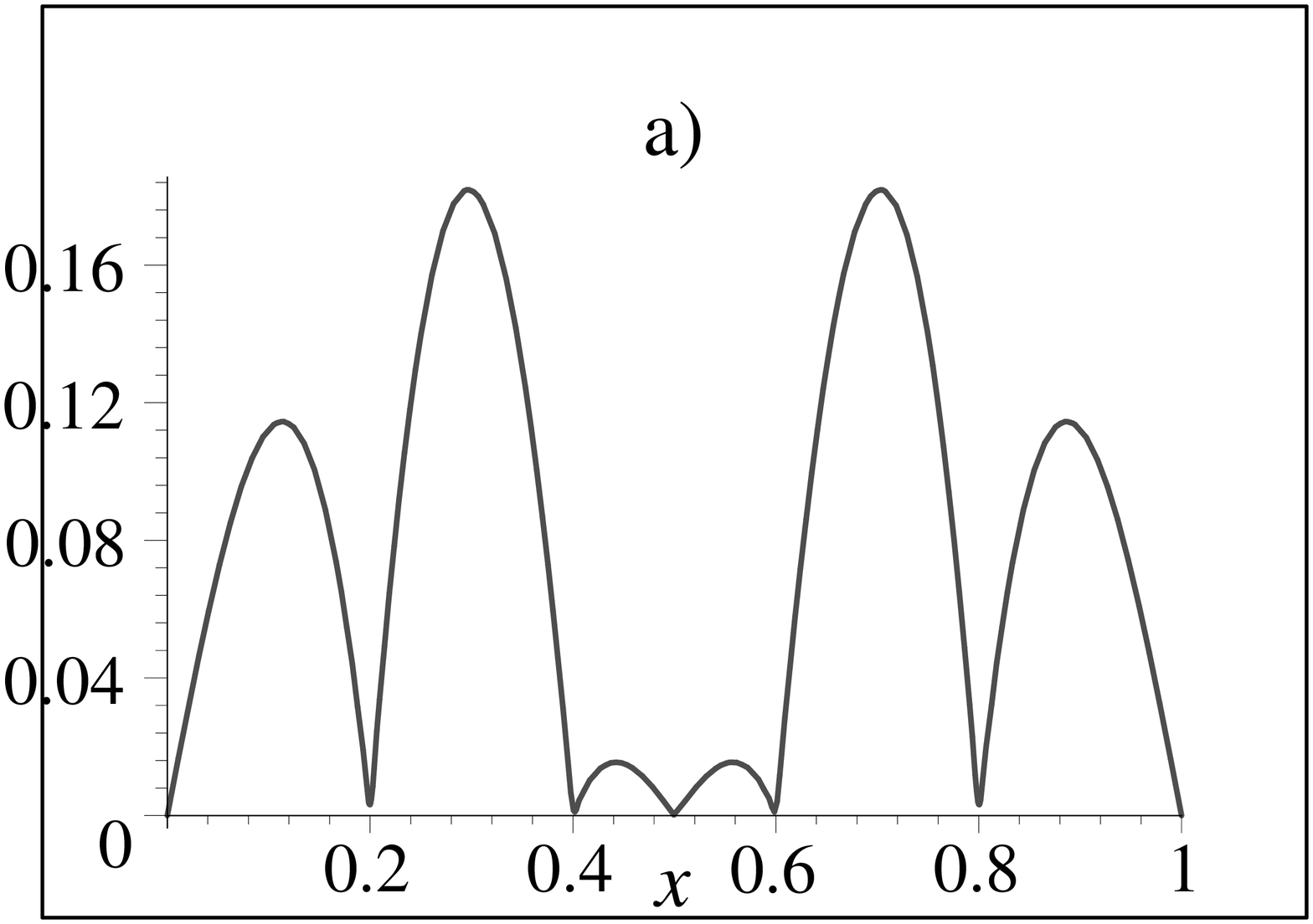}\hspace{1cm}
\includegraphics[width=100pt,angle=-90]{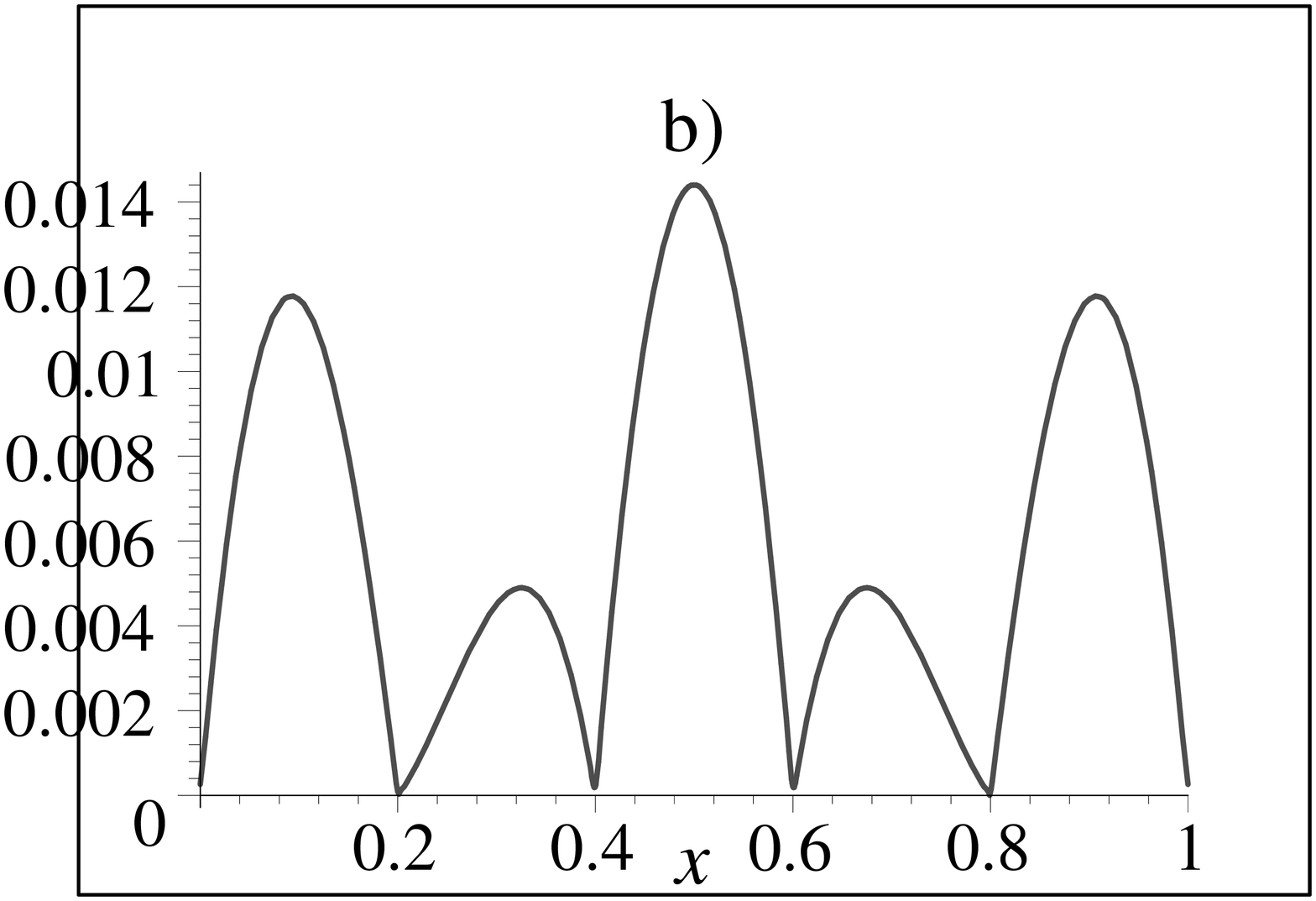}
\caption{Graphs of absolute errors for $m=1$, $N=5$:
 a) $|\sin(2\pi x)-P_{\varphi}(x)|$, b) $|B_{10}(x)-P_{\varphi}(x)|$.} \label{Fig2}
\end{figure}

\emph{The case $m=2$.}
\begin{figure}[h]
\includegraphics[width=60pt,angle=-90]{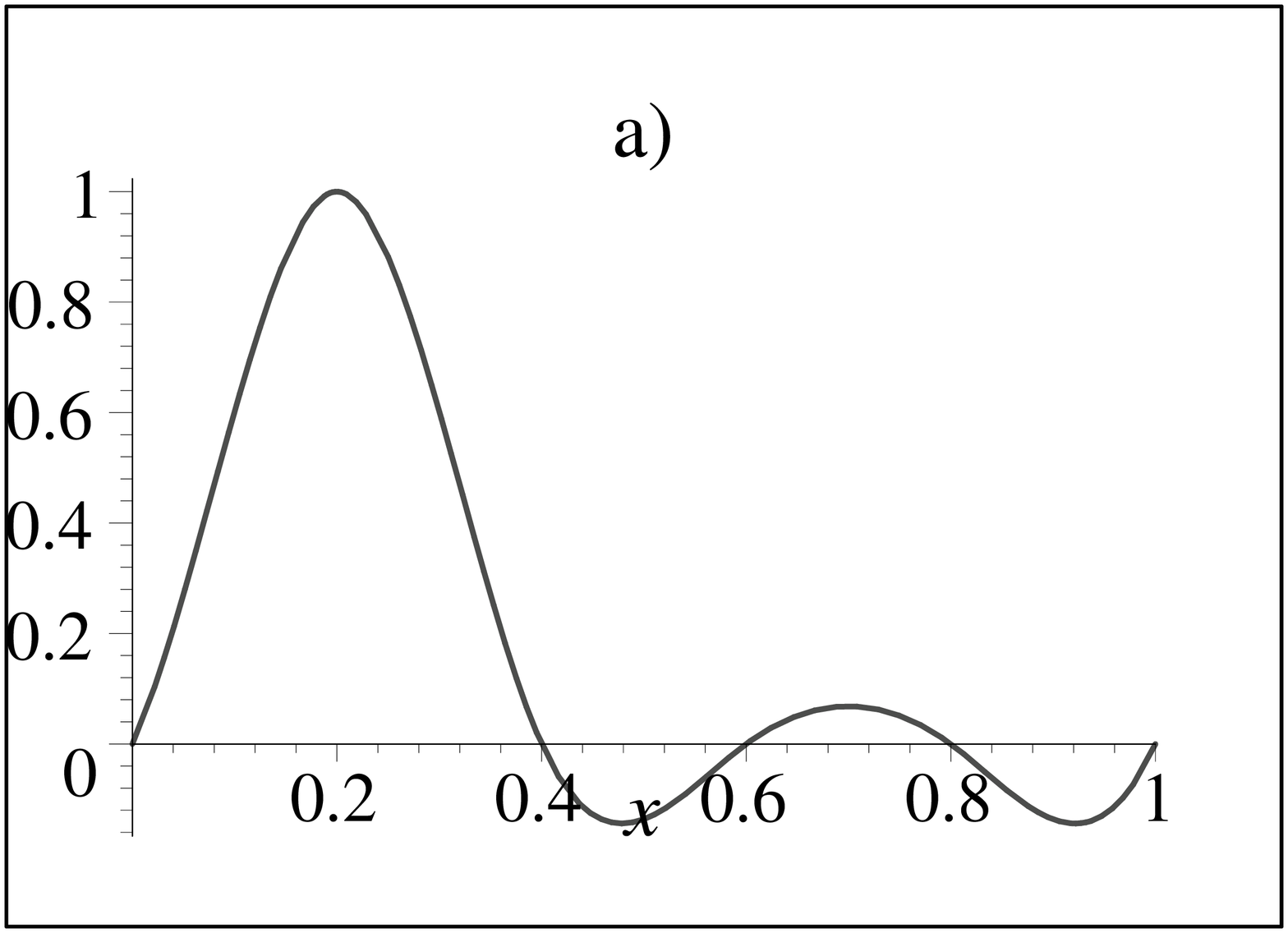}
\includegraphics[width=60pt,angle=-90]{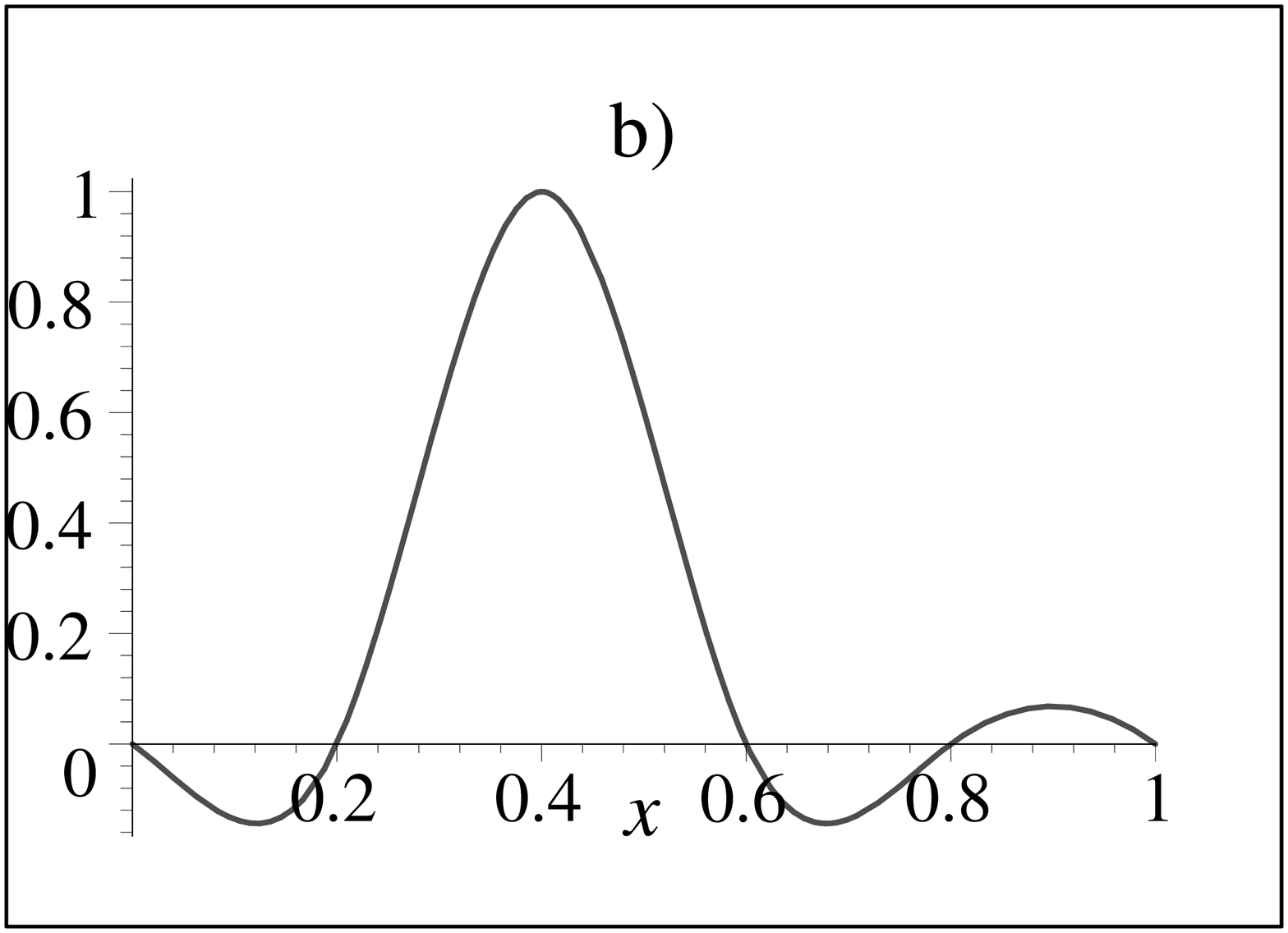}
\includegraphics[width=60pt,angle=-90]{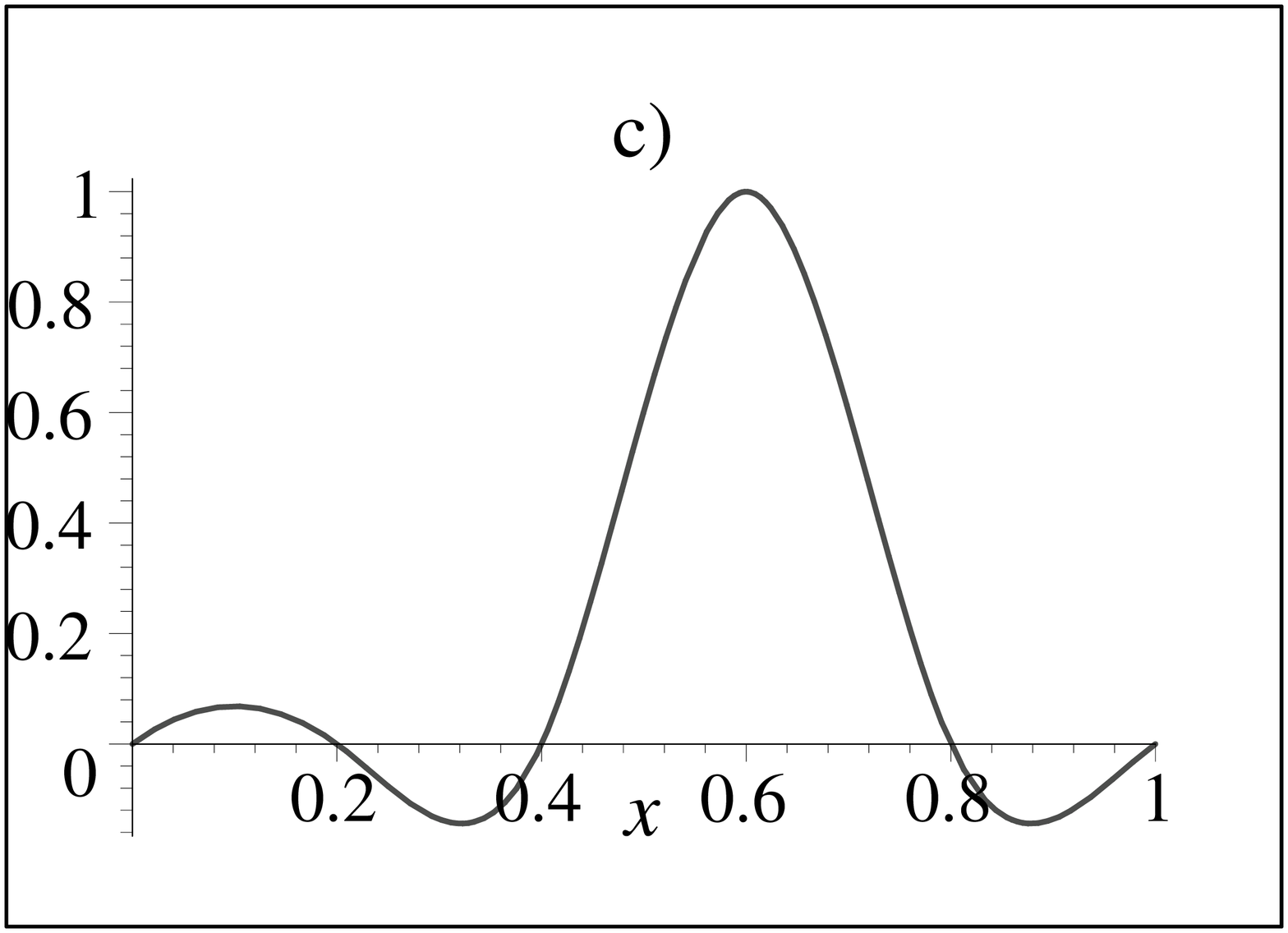}
\includegraphics[width=60pt,angle=-90]{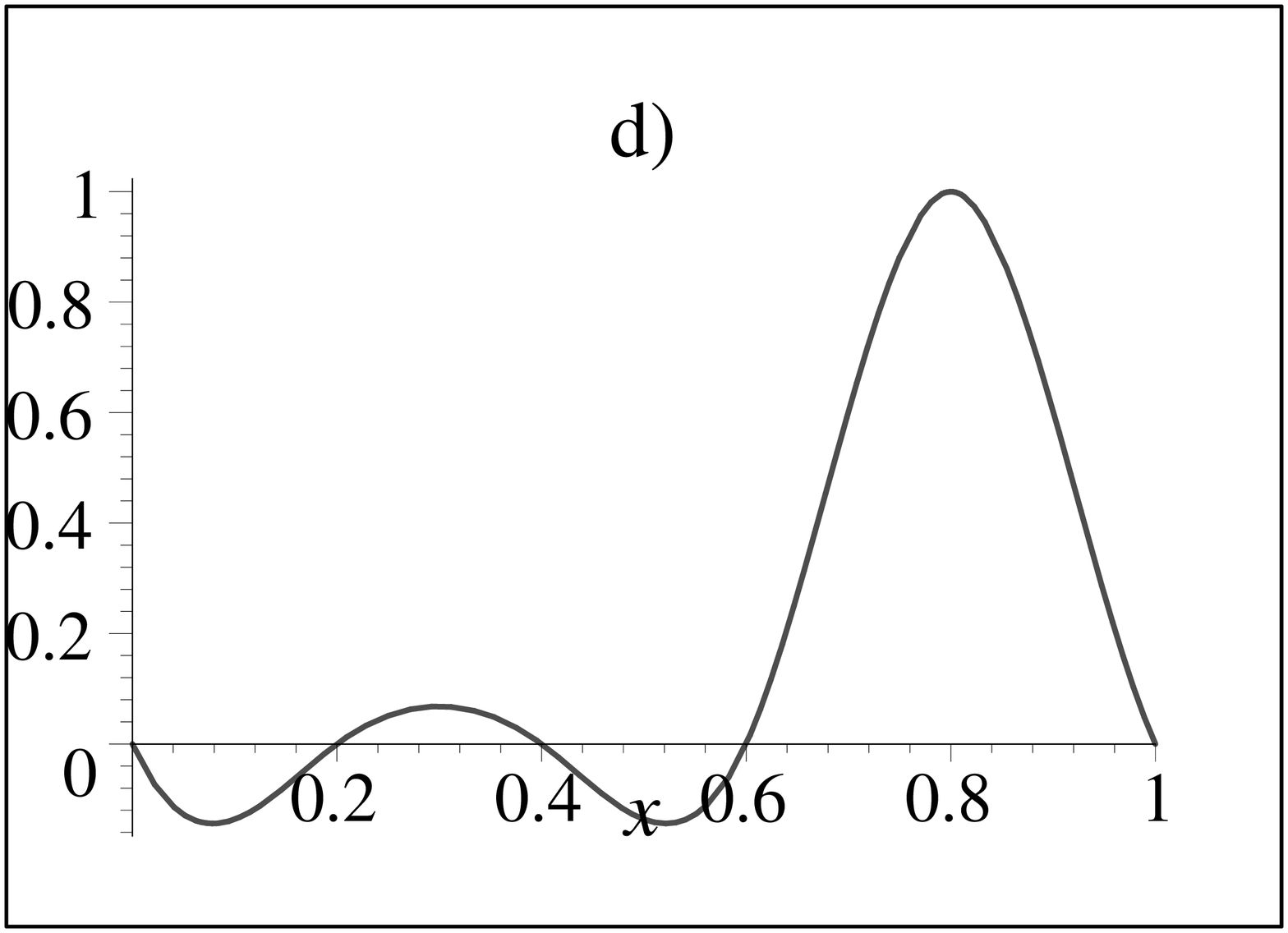}
\includegraphics[width=60pt,angle=-90]{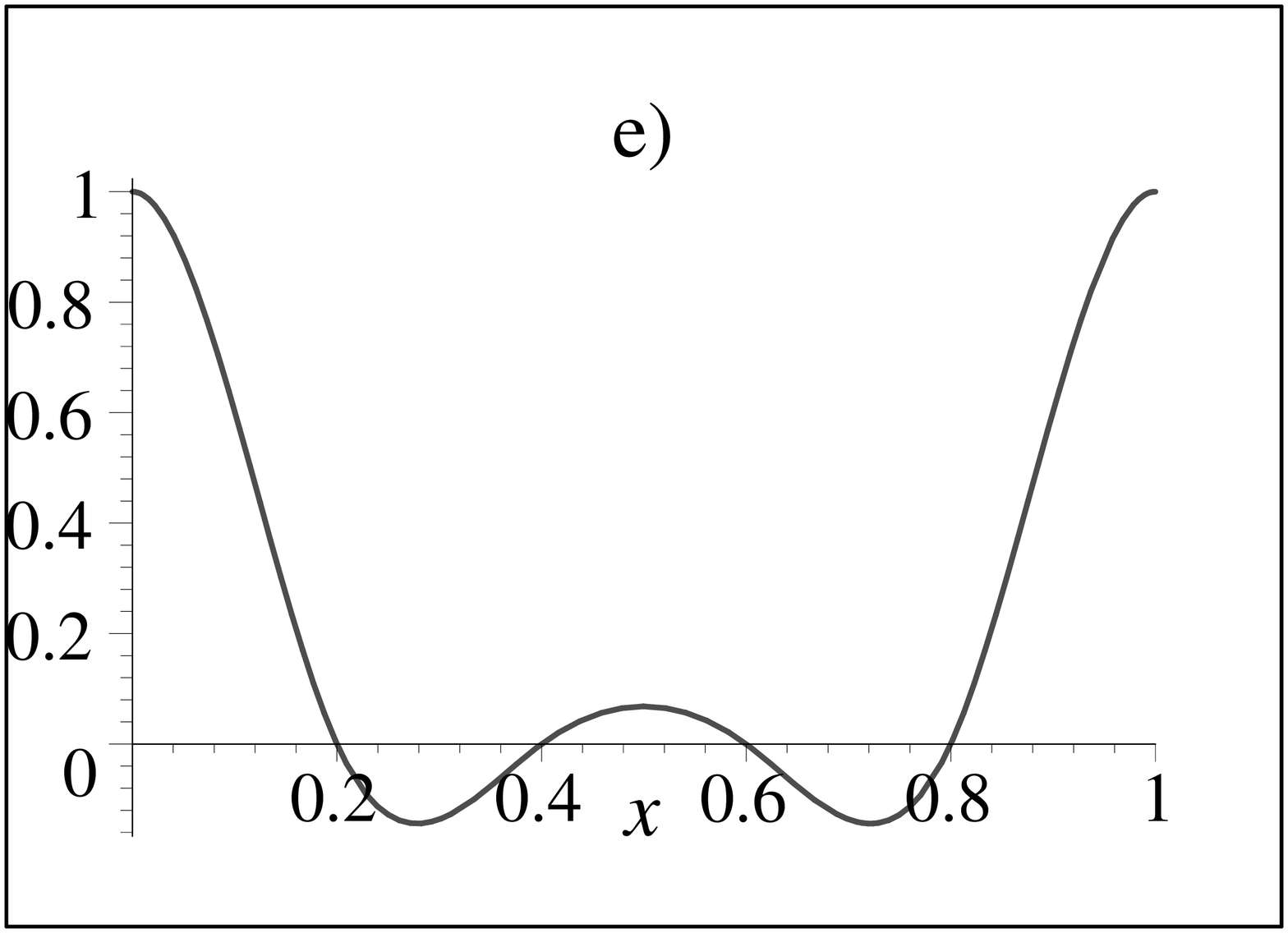}\\
\caption{ Graphs of optimal coefficients for $m=2$, $N=5$:
 a) $C([1],x)$, b) $C([2],x)$, c) $C([3],x)$, d) $C([4],x)$, e)
$C([5],x)$. } \label{Fig3}
\end{figure}
\begin{figure}[h]
\includegraphics[width=100pt,angle=-90]{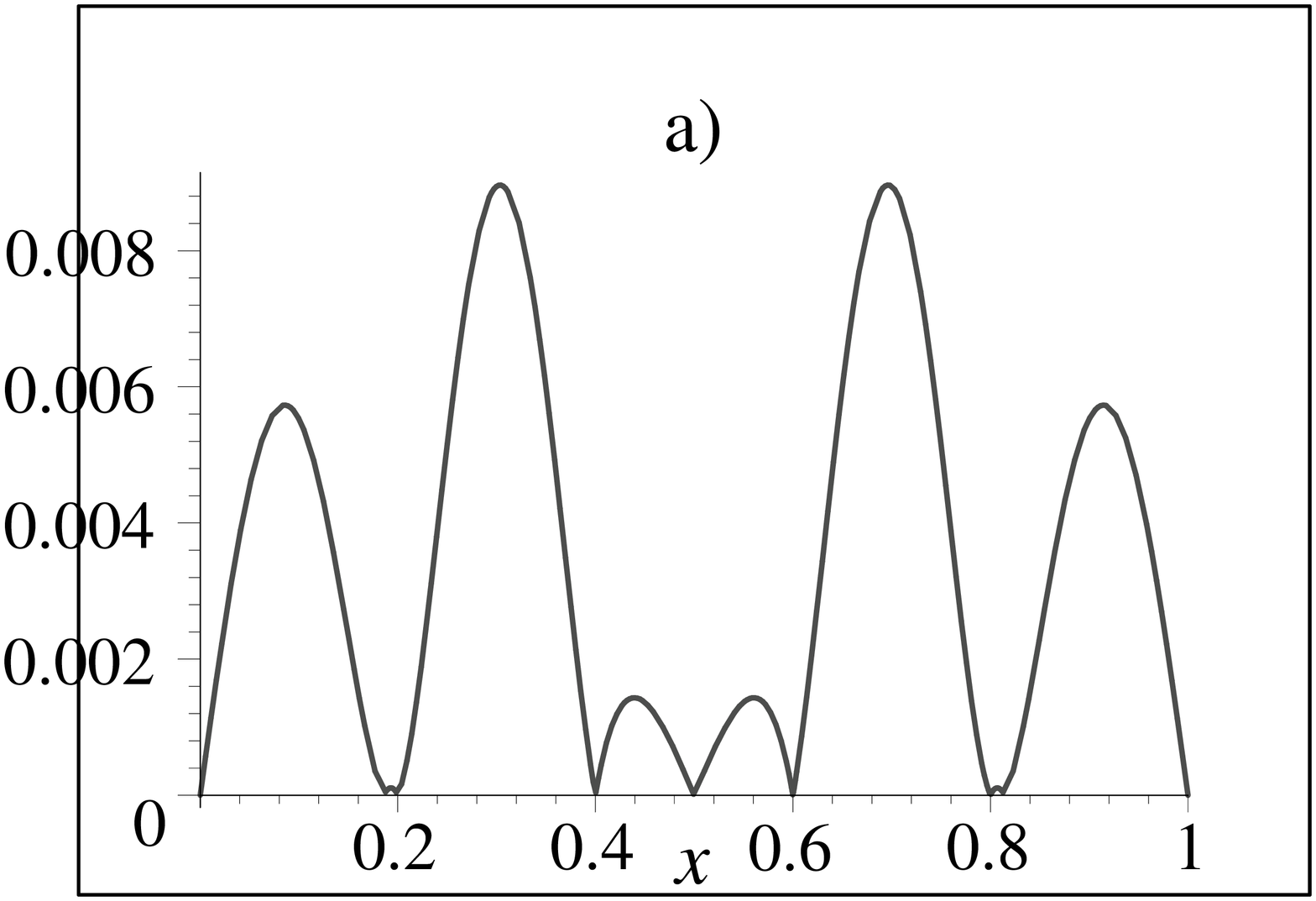}\hspace{1cm}
\includegraphics[width=100pt,angle=-90]{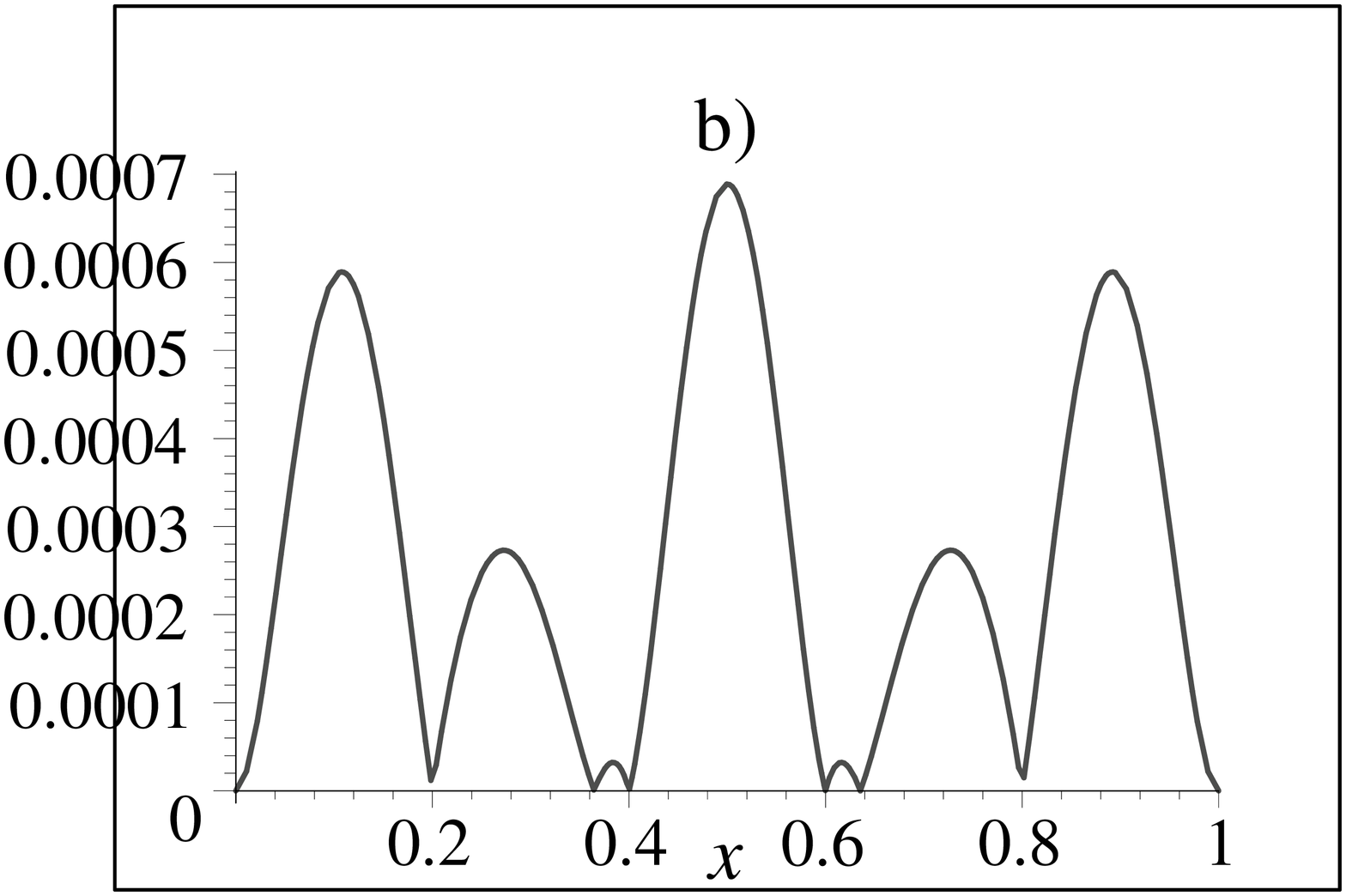}
\caption{Graphs of absolute errors for $m=2$, $N=5$:
 a) $|\sin(2\pi x)-P_{\varphi}(x)|$, b) $|B_{10}(x)-P_{\varphi}(x)|$.} \label{Fig4}
\end{figure}

\emph{The case $m=3$.}
\begin{figure}[h]
\includegraphics[width=60pt,angle=-90]{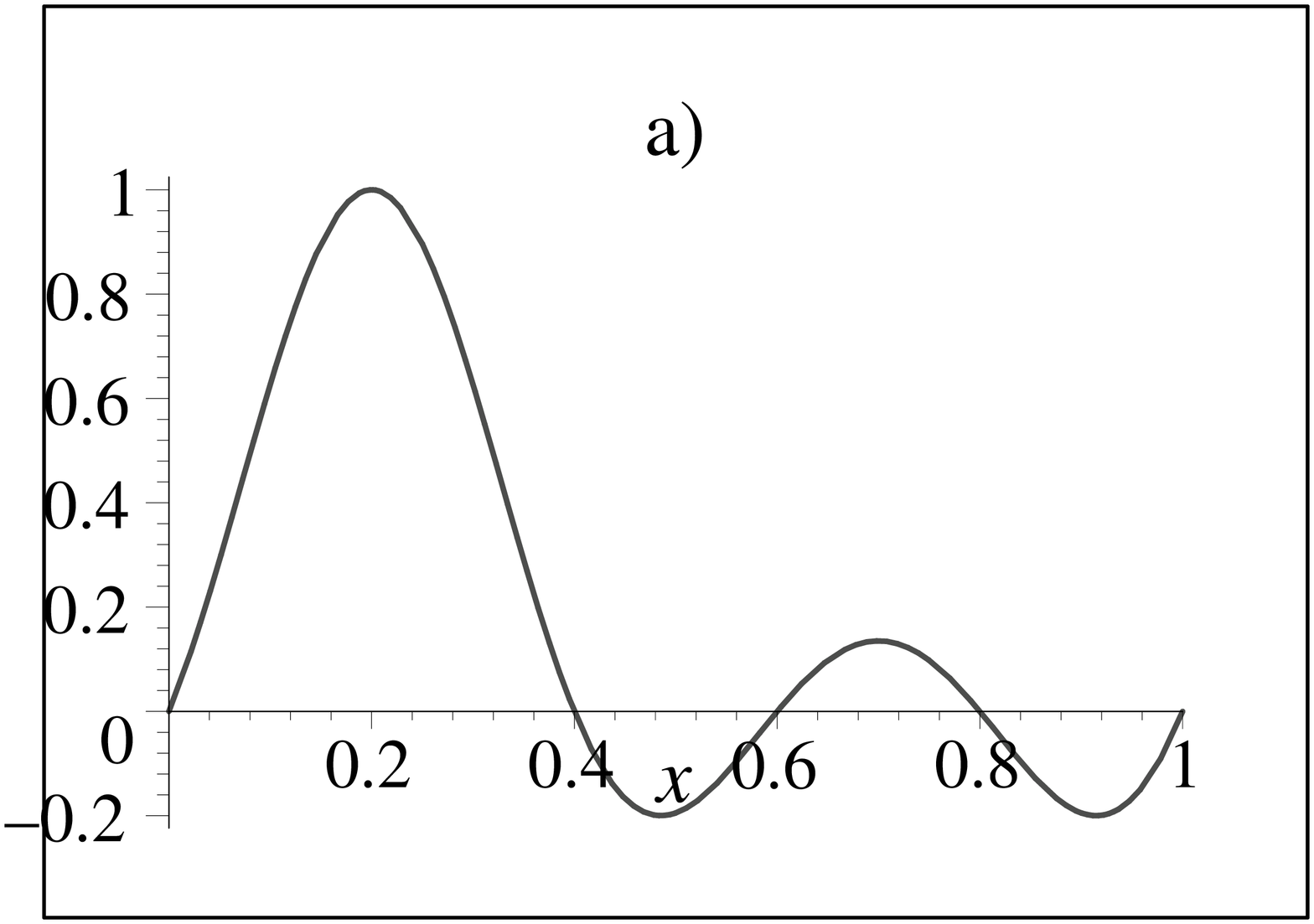}
\includegraphics[width=60pt,angle=-90]{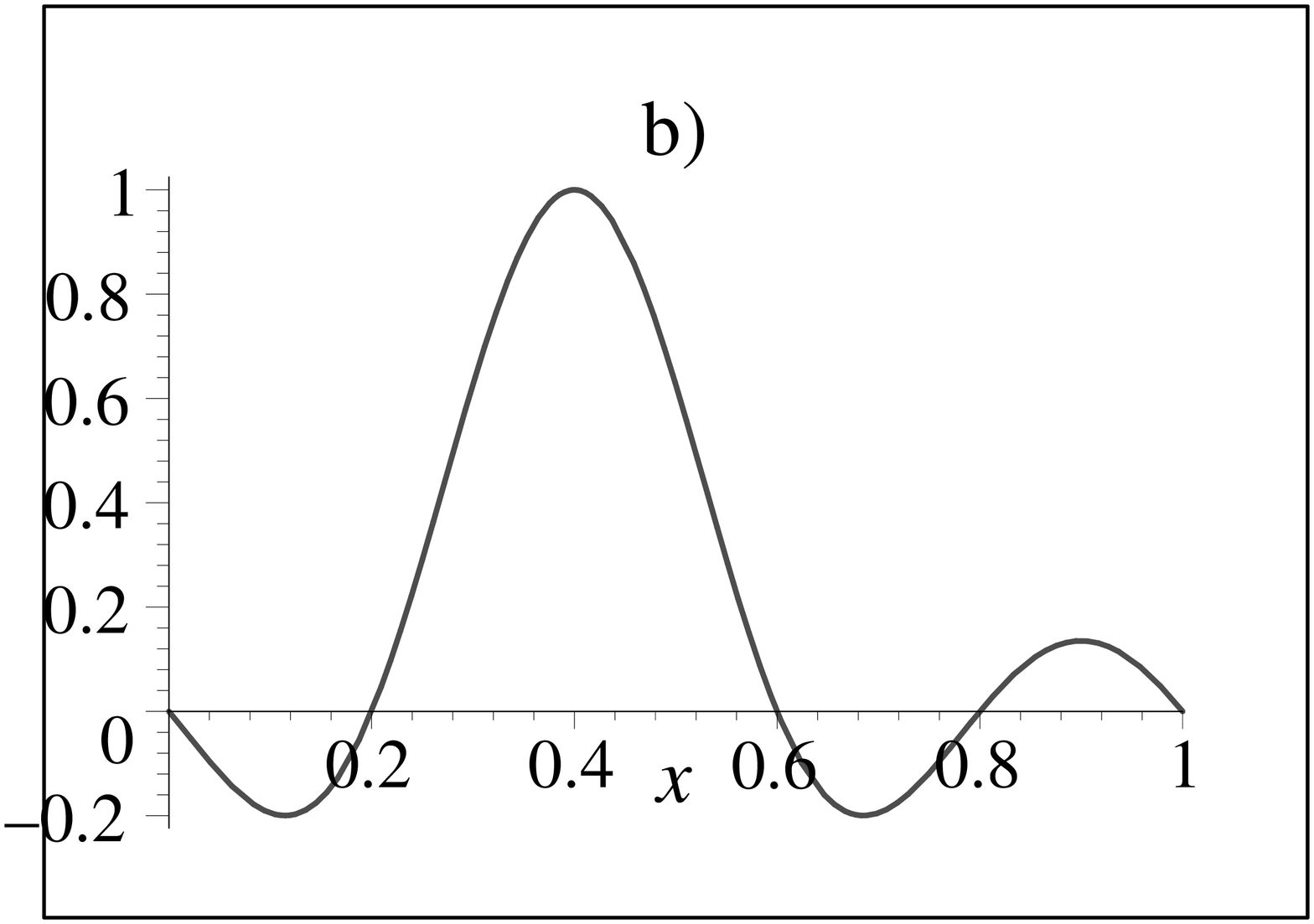}
\includegraphics[width=60pt,angle=-90]{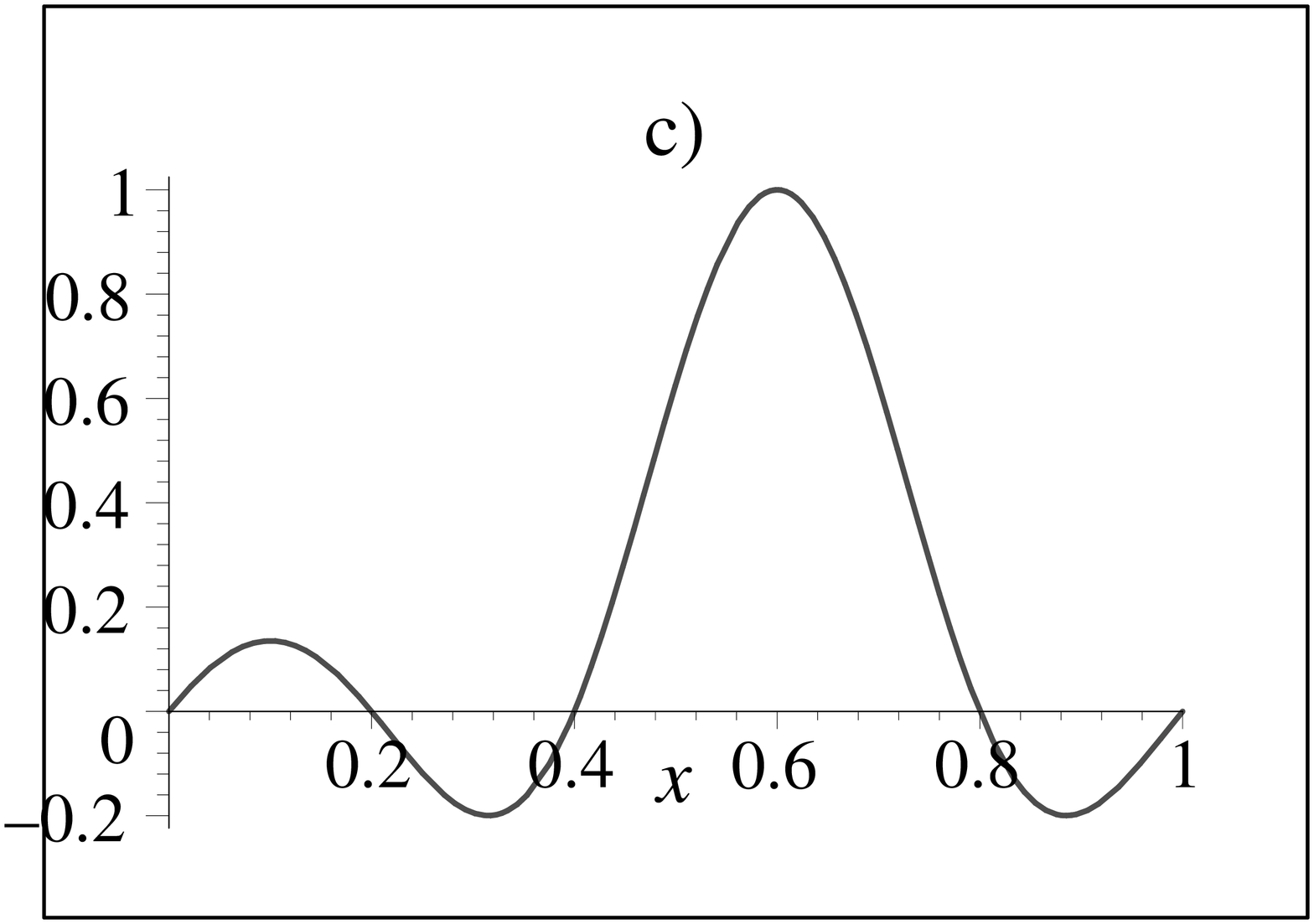}
\includegraphics[width=60pt,angle=-90]{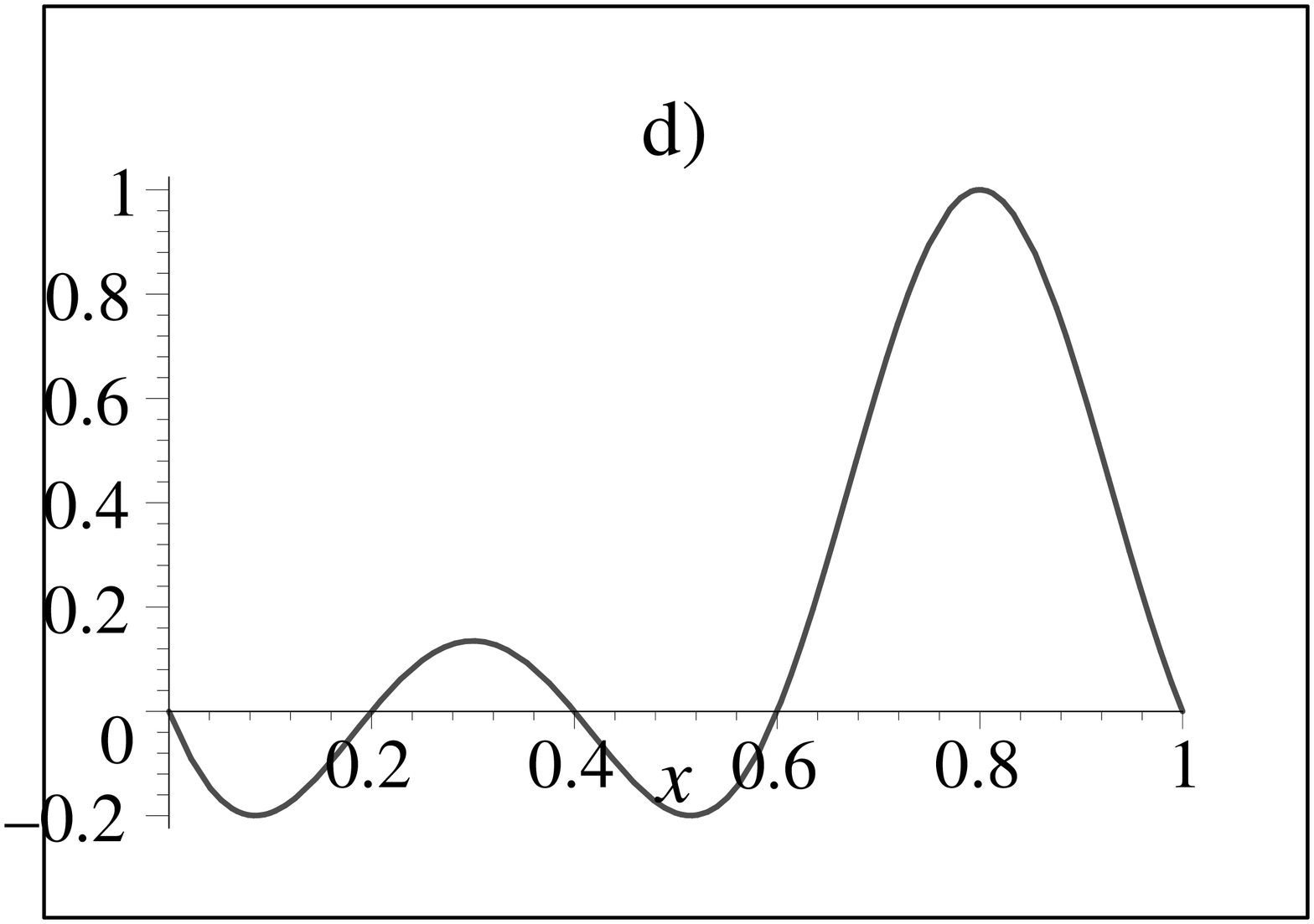}
\includegraphics[width=60pt,angle=-90]{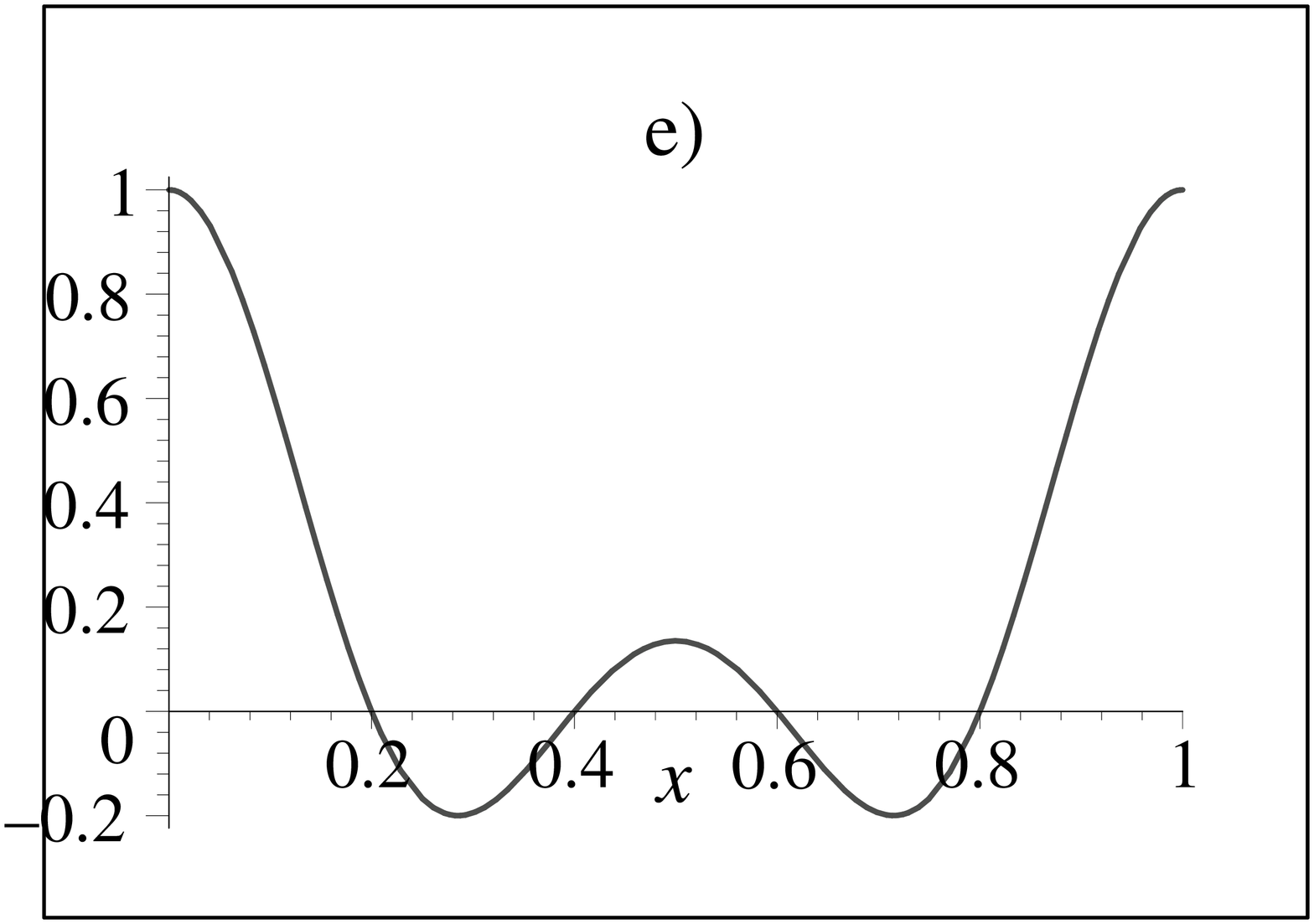}\\
\caption{ Graphs of optimal coefficients for $m=3$, $N=5$:
 a) $C([1],x)$, b) $C([2],x)$, c) $C([3],x)$, d) $C([4],x)$, e)
$C([5],x)$. } \label{Fig5}
\end{figure}
\begin{figure}[h]
\includegraphics[width=100pt,angle=-90]{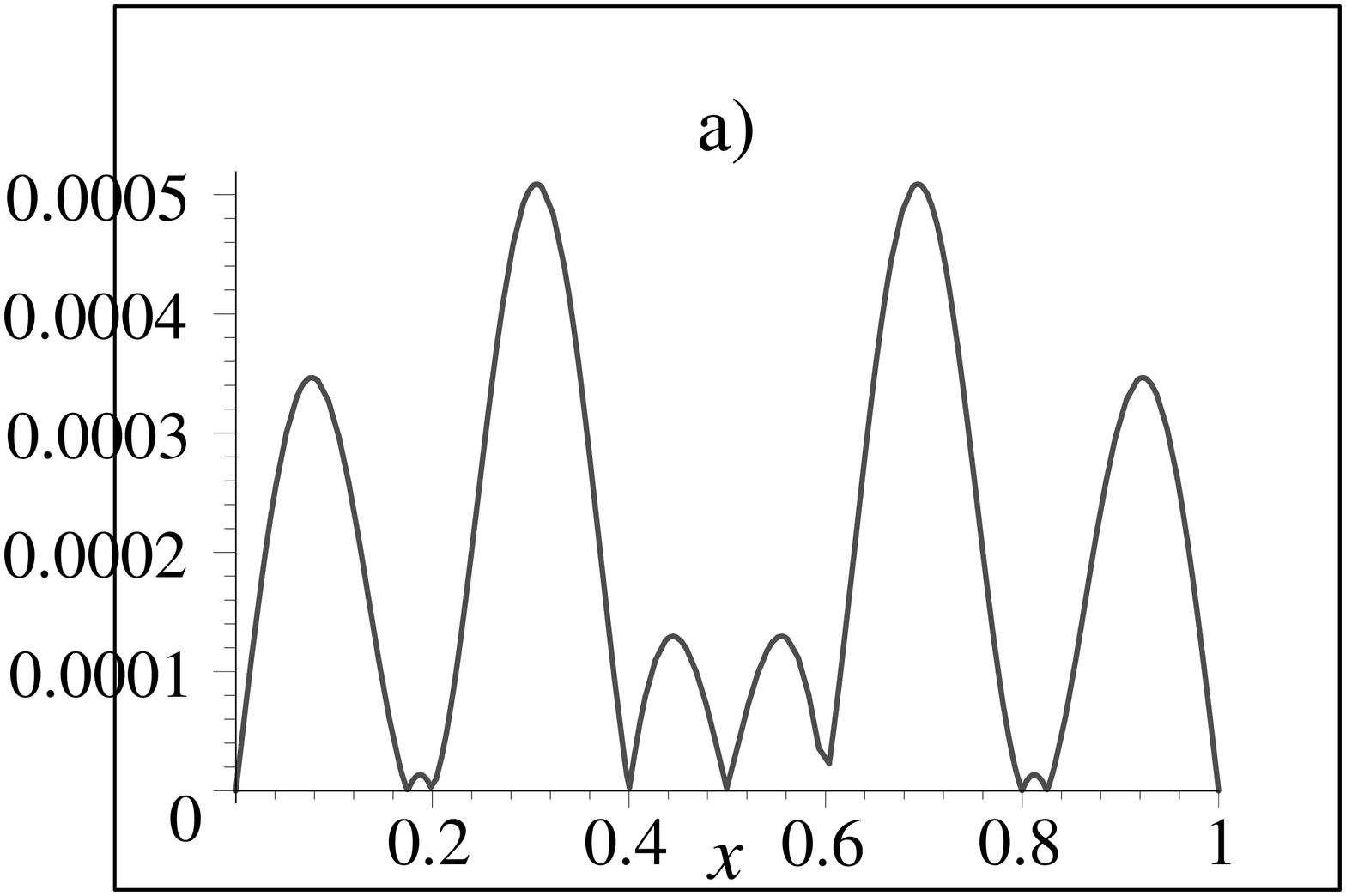}\hspace{1cm}
\includegraphics[width=100pt,angle=-90]{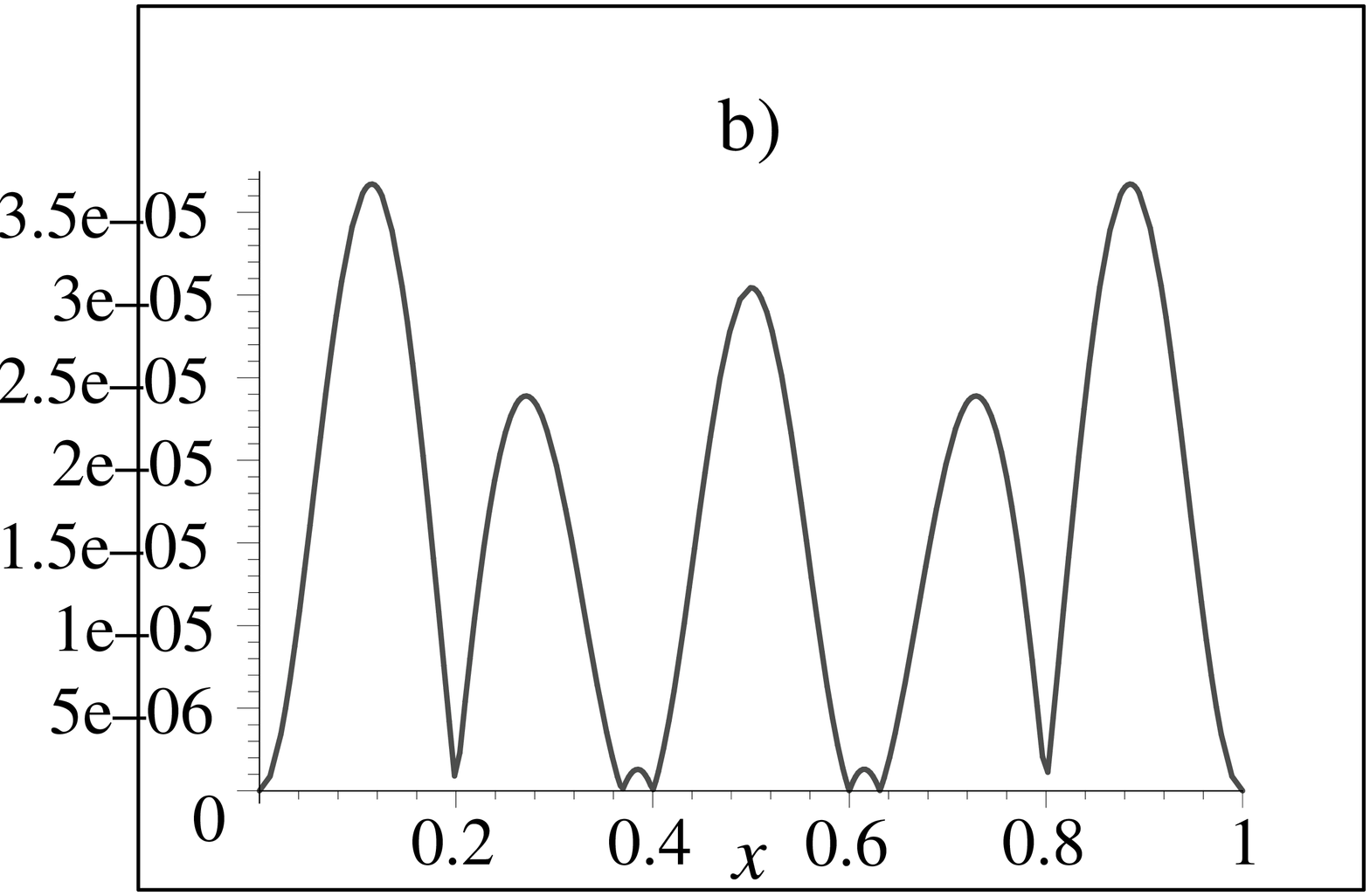}
\caption{Graphs of absolute errors for $m=3$, $N=5$:
 a) $|\sin(2\pi x)-P_{\varphi}(x)|$, b) $|B_{10}(x)-P_{\varphi}(x)|$.} \label{Fig6}
\end{figure}

\emph{The case $m=4$.}
\begin{figure}[h]
\includegraphics[width=60pt,angle=-90]{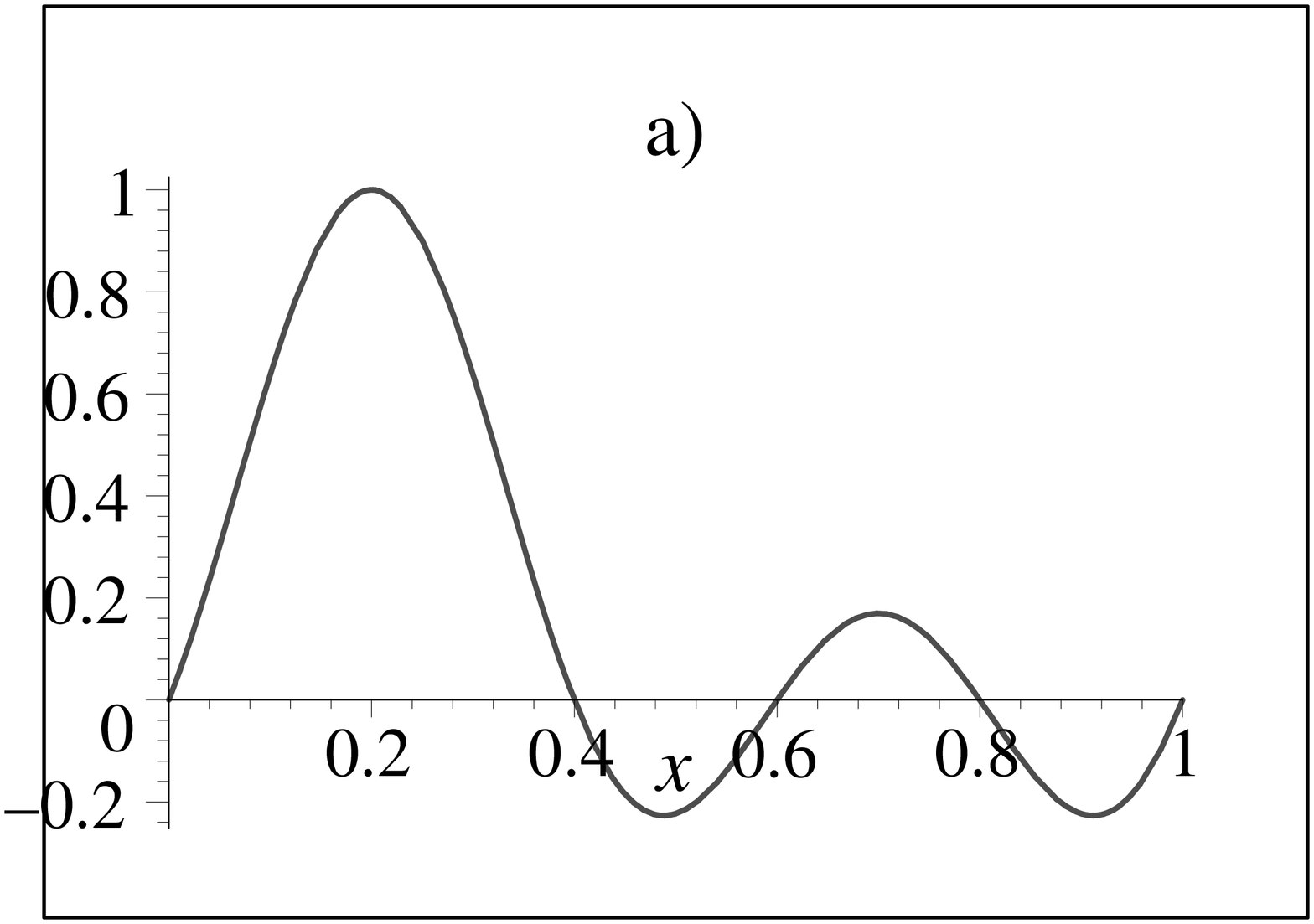}
\includegraphics[width=60pt,angle=-90]{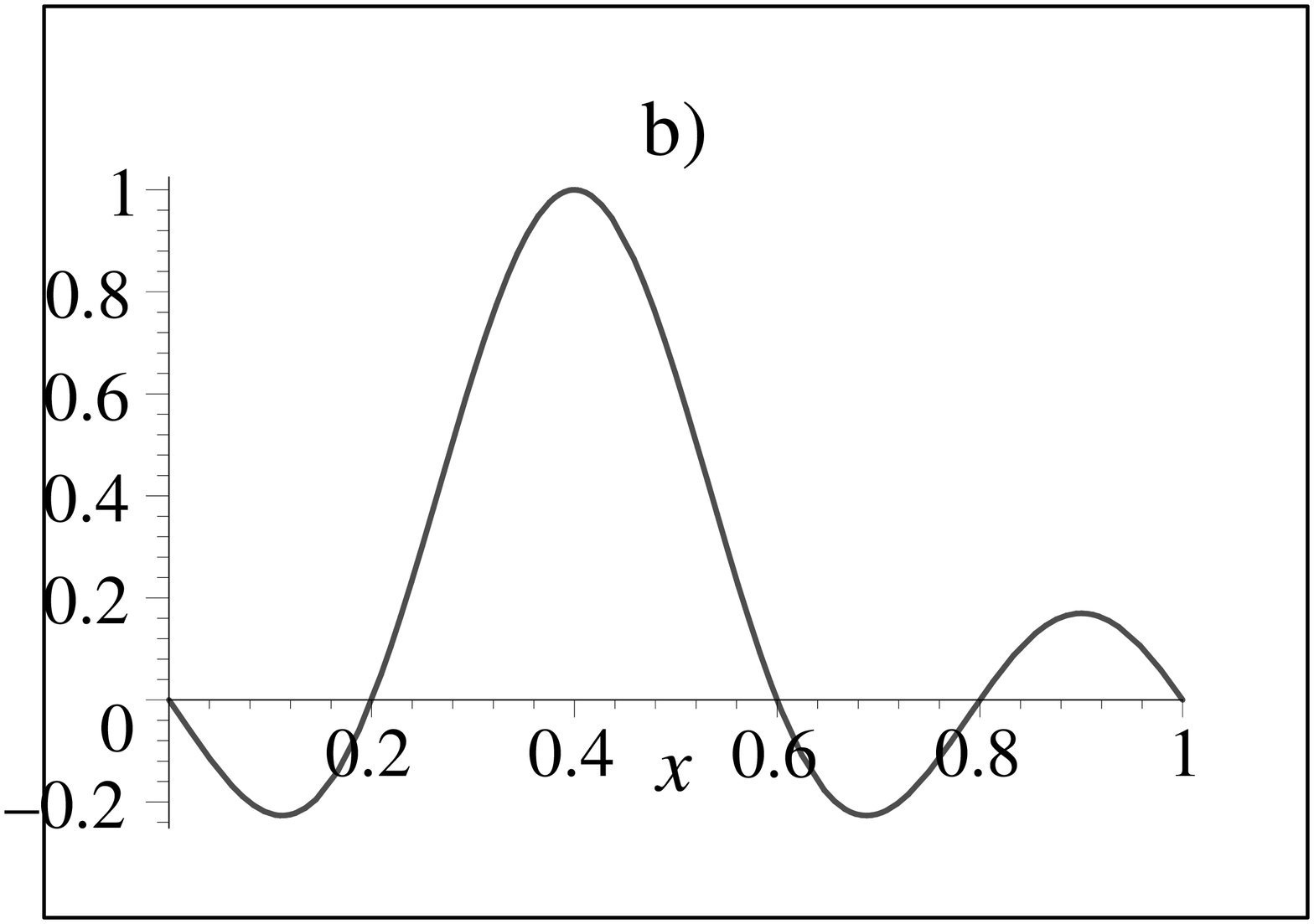}
\includegraphics[width=60pt,angle=-90]{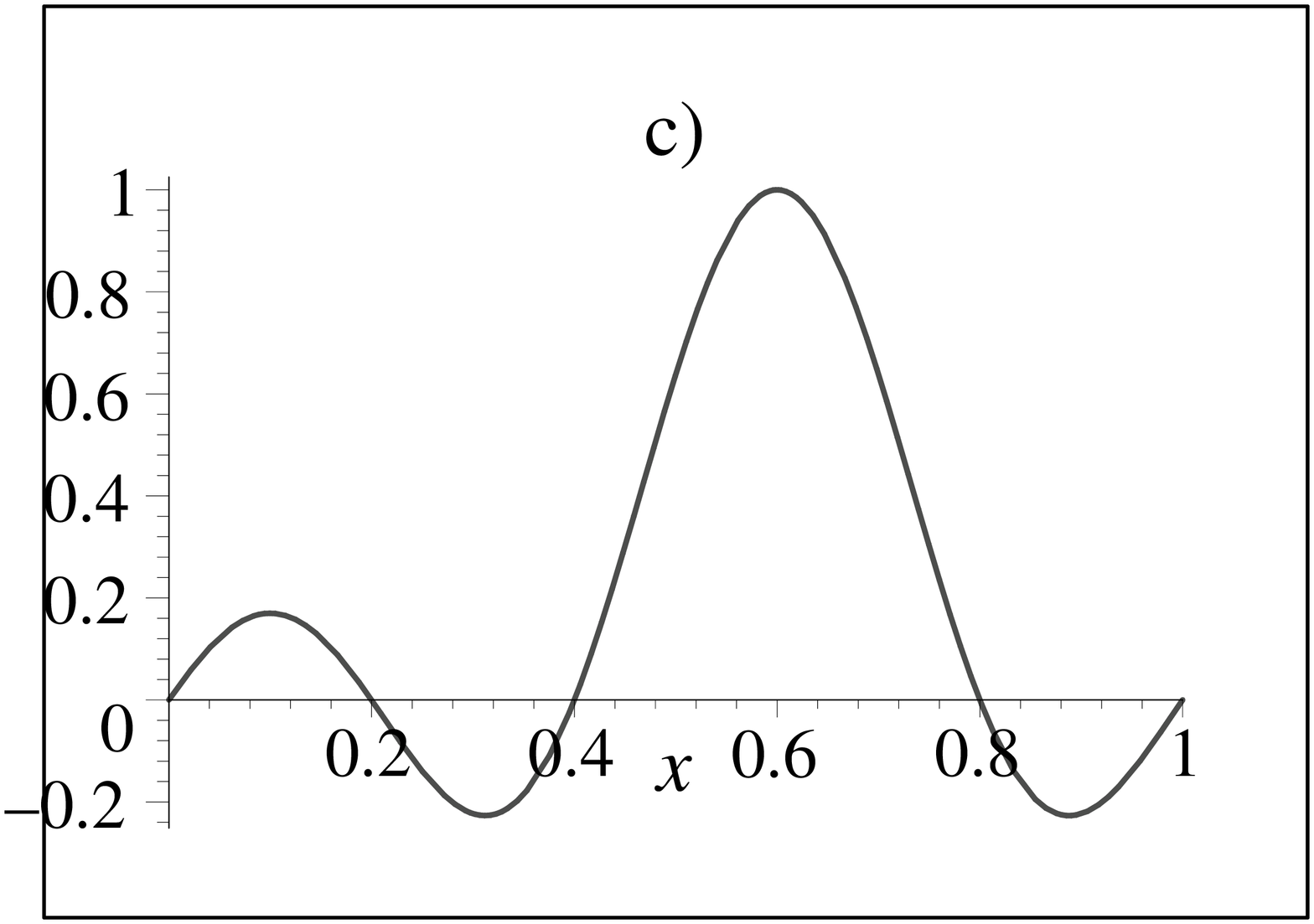}
\includegraphics[width=60pt,angle=-90]{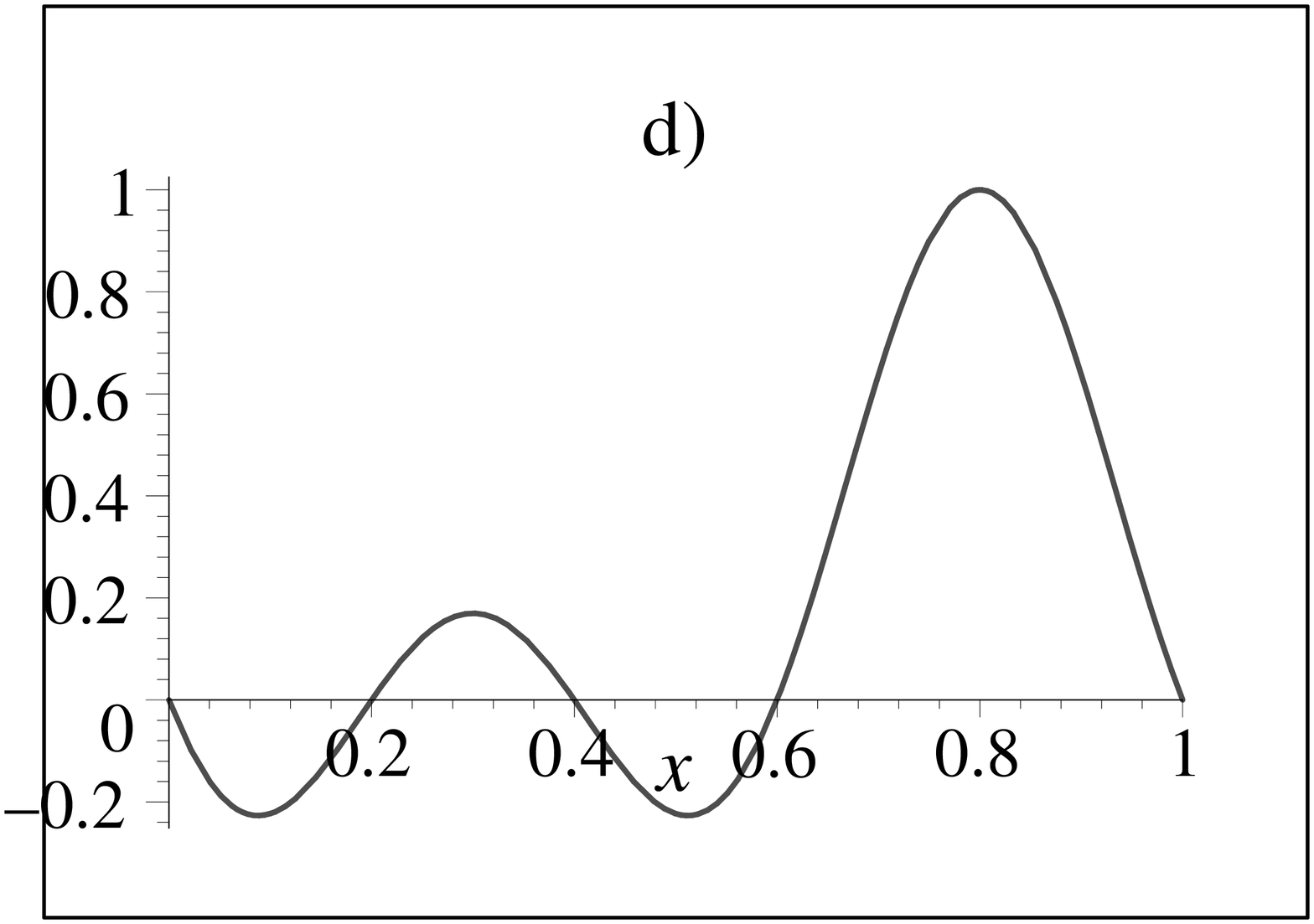}
\includegraphics[width=60pt,angle=-90]{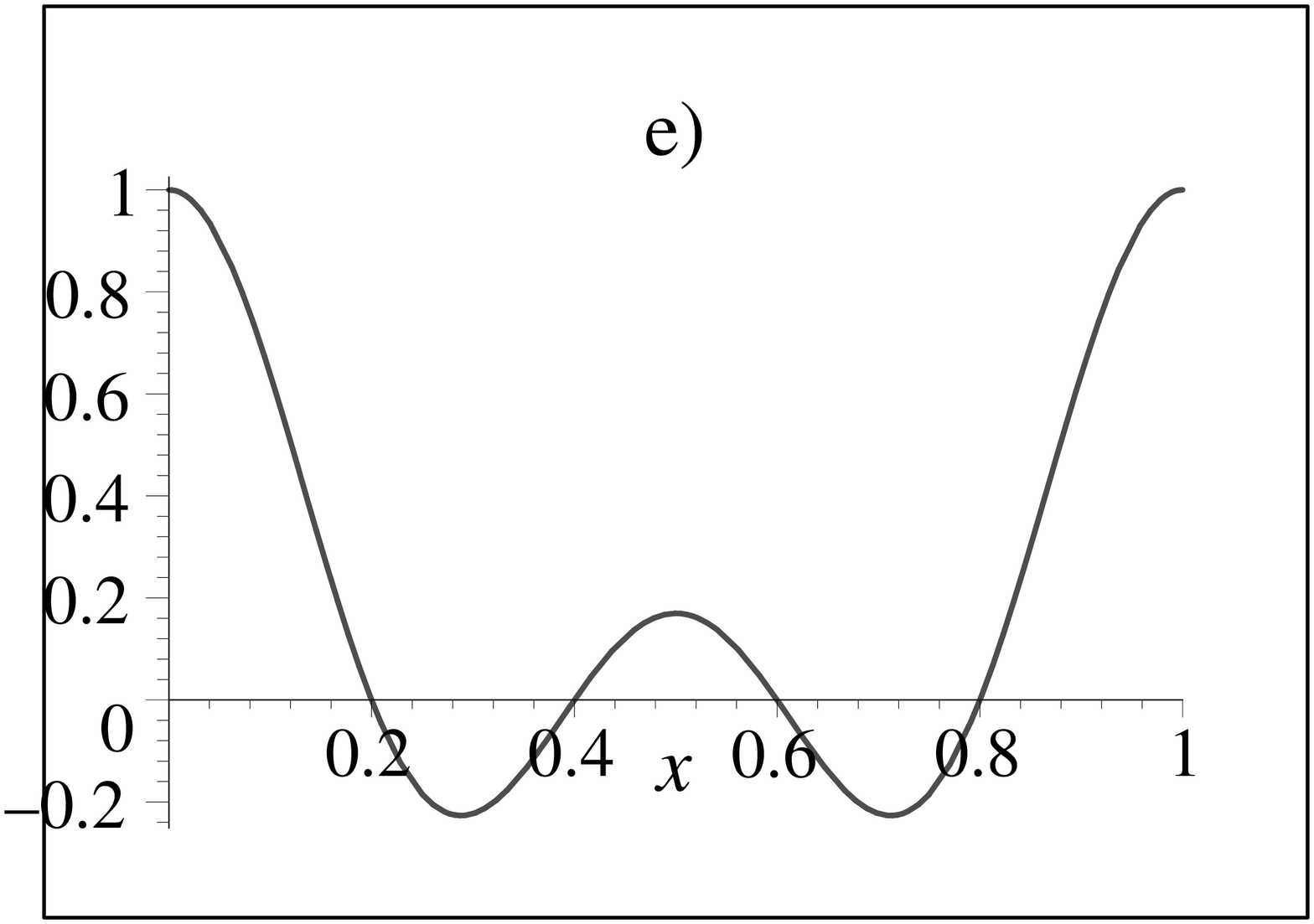}\\
\caption{ Graphs of optimal coefficients for $m=4$, $N=5$:
 a) $C([1],x)$, b) $C([2],x)$, c) $C([3],x)$, d) $C([4],x)$, e)
$C([5],x)$. } \label{Fig7}
\end{figure}
\begin{figure}[h]
\includegraphics[width=100pt,angle=-90]{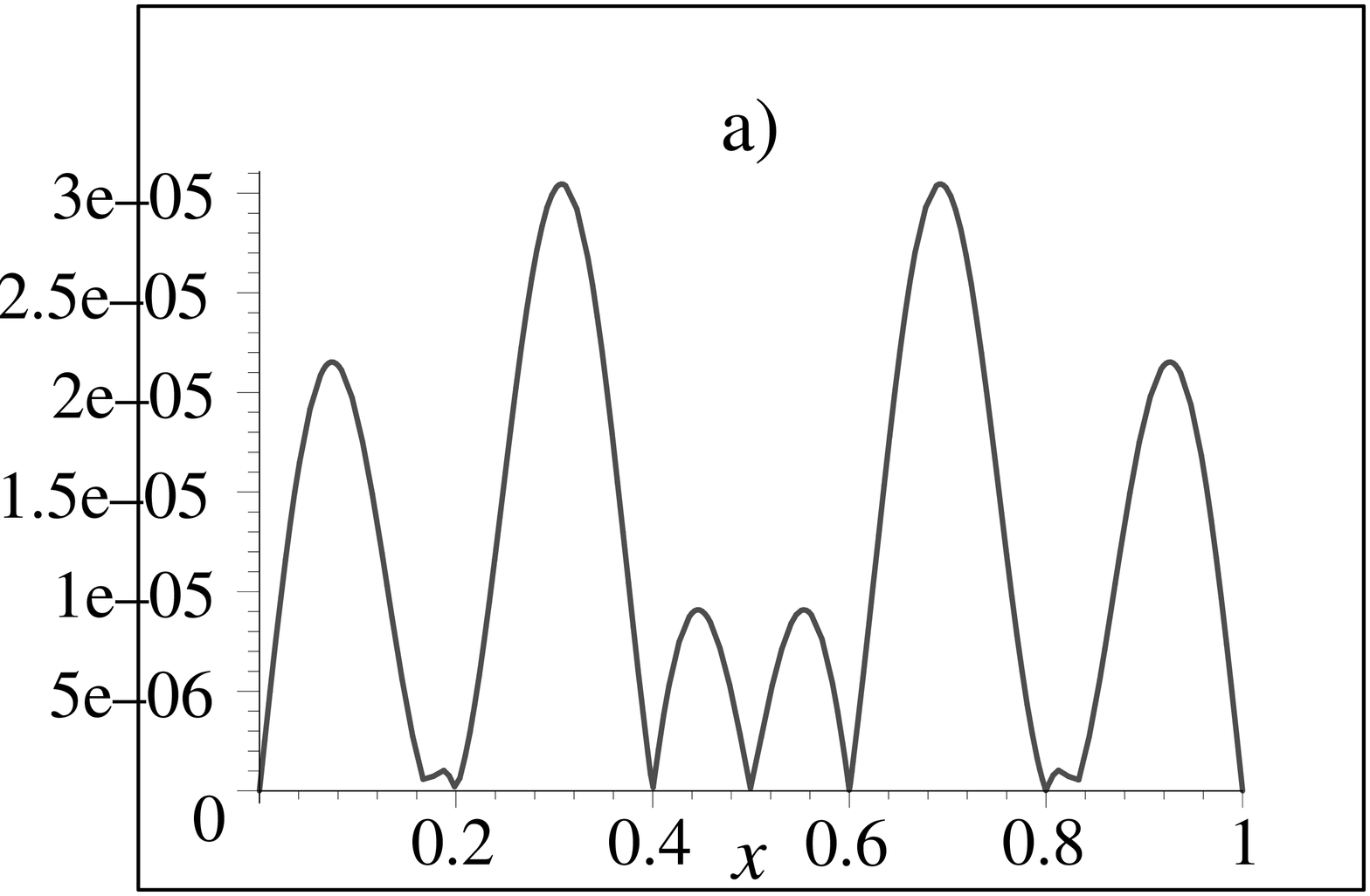}\hspace{1cm}
\includegraphics[width=100pt,angle=-90]{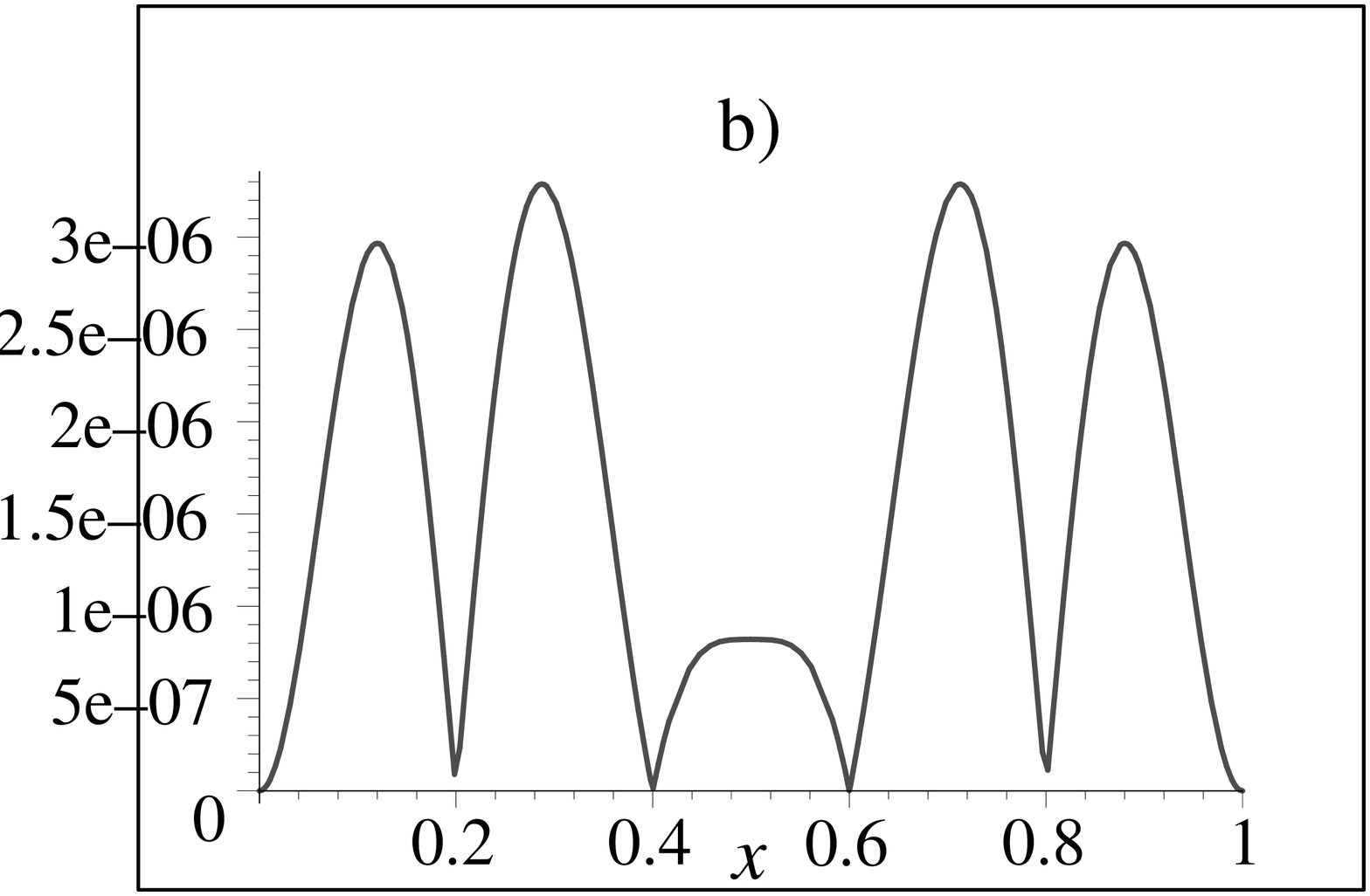}
\caption{Graphs of absolute errors for $m=4$, $N=5$:
 a) $|\sin(2\pi x)-P_{\varphi}(x)|$, b) $|B_{10}(x)-P_{\varphi}(x)|$.} \label{Fig8}
\end{figure}

\newpage
From Figures \ref{Fig2},\ref{Fig4},\ref{Fig6},\ref{Fig8} we can
conclude that the absolute errors between given functions and
optimal interpolation formula is decreasing as $m$ is increasing.

\section*{ Acknowledgements}
The authors are very thankful to professor Erich Novak for
discussion of the results of this paper. A.R. Hayotov thanks
professor Erich Novak and his research group for hospitality.

The final part of this work was done in the Friedrich-Schiller
University of Jena, Germany. The second author thanks the DAAD for
scholarship. Furthermore the second author also thanks
IMU/CDE-program for travel support to the Friedrich-Schiller
University of Jena, Germany.


\begin{thebibliography}{99}
\bibitem{Ahlb67} J.H.Ahlberg, E.N.Nilson, J.L.Walsh, The theory of splines and
their applications, Mathematics in Science and Engineering, New
York: Academic Press, 1967.

\bibitem{Boor78} C. de Boor, A practical guide to splines,
Springer-Verlag, 1978.

\bibitem{Schum81} L.Schumaker, Spline functions: basic theory, J.
Wiley, New-York, 1981.

\bibitem{Lor75} P.-J.Laurent, Approximation and Optimization, Mir, Moscow, 1975,
496 p. (in Russian)

\bibitem{Attea92} M.Attea, Hilbertian kernels and spline functions,
Studies in Computational Matematics 4, C. Brezinski and L.Wuytack
eds, North-Holland, 1992.

\bibitem{Stech76} S.B.Stechkin, Yu.N.Subbotin, Splines in computational mathematics,
Nauka, Moscow, 1976, 248 p. (in Russian)

\bibitem{Vas83} V.A.Vasilenko, Spline-fucntions: Theory, Algorithms,
Programs, Nauka, Novosibirsk, 1983, 216 p. (in Russian)

\bibitem{Arc04} R.Arcangeli, M.C.Lopez de Silanes, J.J.Torrens, Multidimensional
minimizing splines, Kluwer Academic publishers. Boston, 2004, 261
p.


\bibitem{Ign91} M.I.Ignatev, A.B.Pevniy, Natural splines of many variables, Nauka, Leningrad,
1991. (in Russian)

\bibitem{Korn93} N.P.Korneichuk, V.F.Babenko, A.A.Ligun, Extremal
properties of polynomials and splines, Naukovo dumka, Kiev, 1992,
304 p. (in Russian)

\bibitem{Wahba90} G.Wahba, Spline models for observational data.
CBMS 59, SIAM, Philadelphia, 1990.

\bibitem{Eub88} R.L.Eubank, Spline smoothing and nonparametric
regression. Marcel-Dekker, New-York, 1988.

\bibitem{GrSi94} P.J.Green and Silverman, Nonparametric regression
and generalized linear models. A roughness penalty approach.
Chapman and Hall, London, 1994.

\bibitem{BerAgnan04} A.Berlinet and C.Thomas-Agnan, Reproducing
Kernel Hilbert Sapces in Probability and Statistics, Kluwer
Academic Publisher, 2004.

\bibitem{Hol57} J.C.Holladay, Smoothest curve approximation, Math. Tables
Aids Comput. V.11. (1957) 223-243.

\bibitem{deBoor63} C. de Boor, Best approximation propertiesof
spline functions of odd degree, J. Math. Mech. 12, (1963),
pp.747-749.

\bibitem{Schoen64} I.J.Schoenberg, On trigonometric spline
interpolation, J. Math. Mech. 13, (1964), pp.795-825.

\bibitem{Golomb68} M.Golomb, Approximation by periodic spline
interpolants on uniform meshes, Journal of Approximation Theory,
1, (1968), pp. 26-65.

\bibitem{Sob61a} S.L.Sobolev, On Interpolation of Functions of $n$
Variables, in: Selected Works of S.L.Sobolev, Springer, 2006, pp.
451-456.

\bibitem{Sob61b} S.L.Sobolev, Formulas of Mechanical Cubature in
$n$- Dimensional Space, in: Selected Works of S.L.Sobolev,
Springer, 2006, pp. 445-450.

\bibitem{Sob74} S.L.Sobolev, Introduction to the Theory of
Cubature Formulas, Nauka, Moscow, 1974, 808 p.

\bibitem{SobVas} S.L.Sobolev, V.L.Vaskevich. The Theory of Cubature
Formulas. Kluwer Academic Publishers Group, Dordrecht (1997).


\bibitem{SMNikBook88} S.M.Nikolskii, Quadrature Formulas,  Nauka, Moscow,
1988, (in Russian)

\bibitem{Shad85} Kh.M.Shadimetov. The Discrete Analogue of the Differential Operator $d^{2m}
/dx^{2m} $ ant Its Construction,  Vopr. Vychisl. Prikl. Mat.  79,
Tashkent, (1985), 22-35. arXiv:1001.0556 [NA.math]

\bibitem{Sob65} S.L.Sobolev, A Difference Analogue of the Polyharmonic
Equation, in: Selected Works of S.L.Sobolev, Springer, 2006, pp.
529-535.

\end{thebibliography}
\end{document}